\documentclass[11pt,letterpaper]{article}

\def\version{28 December 2022}

\usepackage{cite}
\usepackage[T2A]{fontenc}
\usepackage{inputenc}
\usepackage[english]{babel}
\usepackage{caption}
\inputencoding{cp1251}
\usepackage{morewrites}
\usepackage{etoolbox}

\usepackage{longtable}
\usepackage{array}
\usepackage{comment}

\usepackage[textwidth=165mm,textheight=235mm]{geometry}

\usepackage[normalem]{ulem}

\usepackage{enumitem}
\usepackage{wrapfig}

\usepackage{amsmath,amsthm,amssymb}
\usepackage{graphicx}

\newcommand{\ac}[1]{\noindent\textcolor{red}{{\bf [}$\blacktriangleright\!\!\blacktriangleright$\bf #1]}}
\renewcommand{\ac}[1]{}
\usepackage[usenames,dvipsnames]{color}
\definecolor{MyDarkBlue}{rgb}{0,0.08,0.45}
\usepackage[backref=none]{hyperref}
 \hypersetup{pdfborder={0 0 0},
   colorlinks,
   urlcolor={MyDarkBlue},
   linkcolor={MyDarkBlue},
   citecolor={MyDarkBlue},
   breaklinks=true}
\usepackage{doi}
\providecommand{\eprint}[1]{}
\renewcommand{\eprint}[1]{arXiv:\href{http://arxiv.org/abs/#1}{#1}}

\renewcommand{\Re}{\mathop{\rm Re}}

\newtheorem{lemma}{Lemma}[section]
\newtheorem{theorem}[lemma]{Theorem}
\newtheorem{definition}[lemma]{Definition}
\newtheorem{example}[lemma]{Example}
\newtheorem{remark}[lemma]{Remark}

\newcommand{\cm}{{\rm m}}

\newcommand{\dv}{\mathop{\mathrm{div}}}
\newcommand{\curl}{\mathop{\mathrm{curl}}}

\newcommand*{\mailto}[1]{\href{mailto:#1}{\nolinkurl{#1}}}

\newcommand{\R}{\mathbb{R}}

\providecommand{\C}{\mathbb{C}}
\renewcommand{\C}{\mathbb{C}}

\newcommand{\const}{\mathop{\mathrm{const}}\nolimits}

\numberwithin{equation}{section}

\providecommand{\hsk}[1]{\advance\thickmuskip-#1mu\advance\medmuskip-#1mu}

\makeatletter\@addtoreset{equation}{section}
\makeatletter\@addtoreset{lemma}{section}
\makeatother

\begin{document}

\title{Attractors   of  Hamiltonian nonlinear partial differential equations}

\author{Editors: Andrew Comech, Alexander Komech, Mikhail Vishik
\\~\\
}

\author{
\normalsize\sc Andrew Comech\footnote{This work was supported by a grant from the Simons Foundation (851052, A.C.)}
\\
\small\it Texas A\&M University, College Station, Texas 77843, USA
\\[2ex]
\normalsize\sc Alexander Komech
\\
\small\it Vienna, Austria
\\[2ex]
\normalsize\sc Elena Kopylova\footnote{Supported partly by Austrian Science Fund (FWF):  P27492-N25 (E.K.)}
\\
\small\it Vienna University, Vienna A-1090, Austria
}

\date{\version}

\maketitle


 \bigskip\bigskip

\hfill
{\it To the memory of Mark Vishik}

\bigskip
\bigskip

\begin{abstract}
We survey the theory of attractors of nonlinear Hamiltonian partial differential  equations since its appearance in 1990. These are  results
on global attraction to stationary states, to solitons and to stationary orbits,
on adiabatic effective dynamics of solitons and their asymptotic stability.
Results of numerical simulations are also given.
Based on these results,
we propose
a new general hypothesis on attractors
of $G$-invariant nonlinear Hamiltonian partial differential equations.
The obtained results  suggest a novel dynamical interpretation 
of basic quantum phenomena: Bohr's transitions between quantum stationary states, wave-particle duality, and probabilistic interpretation.
\end{abstract}
 


\bigskip

\section{Introduction} 

Theory of attractors of nonlinear PDEs originated from the seminal paper of Landau \cite{L1944} published in 1944, where he
suggested the first mathematical interpretation of turbulence as the growth of the dimension of attractors of the Navier--Stokes equations 
when the Reynolds number increases.
Modern development of the theory of attractors for general \textit{dissipative systems}, i.e., systems with friction 
(the Navier--\allowbreak Stokes equations,  nonlinear parabolic equations, reaction-diffusion equations, wave equations with friction, etc.)
originated in the 1975--1985's in the works of Foias, Hale, Henry, Temam, and others \cite{FMRT2001}--\cite{T1997}.
M.I. Vishik together with his collaborators
worked on global attractors of dissipative nonlinear PDEs from 1980 until 2012. 
These investigations had significantly advanced this theory \cite{BV1992}--\cite{CPV2008}.

A typical result of this theory, in the absence of external excitation, is a global convergence to stationary states:
for any finite energy solution to  dissipative {\it autonomous} equation in a~region $\Omega\subset\R^n$, there is a convergence
\begin{equation}\label{at11}
	\psi(x,t) \to S(x), \qquad t\to + \infty.
\end{equation}
Here $S(x)$ is a stationary  solution, depending on the initial state $\psi(x,0)$.
As a rule, the convergence holds in the $L^2(\Omega)$-metric.
In particular, the relaxation to an equilibrium regime in chemical reactions is due to the energy dissipation.
\medskip

A development of a similar theory for  \textit{Hamiltonian PDEs}
looked unmotivated and even impossible since the convergence (\ref{at11})
seemed to contradict the energy conservation 
and time reversibility admitted by these equations.  On the other hand,   
such a
theory of global attractors is indispensable  for mathematical 
foundation of quantum physics, since all fundamental PDEs  of quantum theory are Hamiltonian: the Maxwell,  Schr\"odinger, Klein--Gordon, Dirac, Yang--Mills equations, etc. 
Even the shape of the theory is already suggested
by fundamental postulates of quantum physics. In particular, the global attraction of type
 (\ref{at11}) is suggested by Bohr's postulates of 1913 on transitions between quantum stationary states.
 Namely, the postulate  can be interpreted as the global attraction 
of all quantum trajectories to an attractor formed by quantum stationary states.
Such a global convergence also allows us to give a mathematical 
interpretation to the 
 L. de Broglie wave particle duality and to
the M. Born probabilistic interpretation; for more details, see
\cite{K2016,K2013,K2021,K2021qj}.
We note that in 1961 W. Heisenberg began developing 
a nonlinear theory of elementary particles \cite{Heis1961,Heis1966}.
\medskip

Thus, the basic postulates of quantum physics
suggest the type of long-time behavior of solutions of
the related  fundamental dynamical equations: Maxwell--Schr\"odinger,
Maxwell--Dirac, Yang--Mills--Dirac, and other equations and systems.
All these coupled equations
are nonlinear Hamiltonian systems of partial differential equations,
and it seems unlikely
that these
and other fundamental dynamical equations
of quantum physics are exceptional among generic Hamiltonian PDEs.
This situation
leads to the conjecture
that the global convergence to a proper attractor
is an inherent feature
of general nonlinear Hamiltonian PDEs.
\smallskip 

The first results on attractors for the Hamiltonian equations  were obtained by  Morawetz,
Segal, and Strauss  in the case when the attractor 
consists of the single point zero. In this case the attraction is called {\it local energy decay}
\cite{Mor1968}--\cite{St81-2}.

It was not until 1990
when
the case of a nontrivial attractor was considered,
when the hypothesis on the global attractors
was set forth by one of the authors.
Since then, this hypothesis
has been developed 
in collaboration with the international  group of researchers:
with M. Kunze, H. Spohn, and B. Vainberg
and later with V.S. Buslaev,
A. Comech,  V. Imaykin, E. Kopylova,
D. Stuart, and others.

The investigations
since 1990 suggest that such long-time behavior of solutions
is not merely the peculiarity of the dynamical equations of quantum physics,
but rather a characteristic property   
of {\it generic} 
nonlinear Hamiltonian PDEs; see the surveys \cite{K2016,KK2020,KK2021}.
Below, we sketch these results. All these results were discussed in many talks
by the authors at  the seminars of M.I. Vishik in the Department of Mechanics and Mathematics
of  Moscow State University.

This theory differs significantly from the theory of attractors of dissipative systems
where the attraction to stationary states is due to an energy dissipation caused by a friction.
For Hamiltonian equations the friction and energy dissipation are absent, 
and the attraction is caused by radiation which irrevocably  brings the energy to infinity. 
It is exactly this radiation that is addressed 
 in the Bohr postulate on transitions between quantum stationary states.

\section{Bohr's postulates:  quantum jumps}
   In 1913, Bohr formulated the following two fundamental postulates of the quantum theory of atoms \cite{Bohr1913-ak}:
\smallskip

\noindent
\ {\bf I.} {\it Each electron lives on one of the quantum stationary orbits, and
sometimes it jumps from one stationary orbit to another:
in the Dirac notation,}
\begin{equation}\label{B1}
| E_n \rangle \mapsto | E_{n'} \rangle.
\end{equation}

\noindent
{\bf II.} 
{\it
 On a stationary orbit\index{stationary orbit},
 the electron does not radiate,
while every  jump
is followed by the radiation of an electromagnetic wave
of frequency 
\begin{equation}\label{B21}
\omega_{nn'} = \frac {E_{n'}-E_n} \hbar = \omega_{n'}-\omega_n, \qquad \omega_n: = E_n / \hbar.
\end{equation}
}
Both these postulates were inspired by stability of atoms, by the Rydberg--Ritz {\it Combination Principle}, and by the Einstein theory of the photoeffect.

In 1926,
with the discovery of the Schr\"odinger theory \cite{S1926},
the question arose about the implementation of the above Bohr axioms in the new theory
based on partial differential equations.
This and other questions
have been frequently addressed in the 1920s and 1930s in heated discussions
by
 N. Bohr\index{Bohr}, E. Schr{\"o}dinger\index{Schr{\"o}dinger!}, A. Einstein\index{Einstein}, and others \cite{B1949}.
However, a satisfactory solution was not achieved, and
a rigorous dynamical interpretation of these postulates is still
unknown.
This lack of theoretical clarity hinders a further progress in the theory (e.g., in superconductivity and in nuclear reactions\index{nuclear! reaction}) and in
numerical simulation\index{numerical simulation} of many engineering processes (e.g., laser radiation\index{laser radiation} and quantum amplifiers\index{quantum! amplifier}), since a computer can solve dynamical equations, but cannot take into account postulates.

\section{Schr\"odinger's identification of stationary orbits}
Besides the equation for the wave function, the Schr\"odinger theory contains a highly
nontrivial definition of stationary orbits 
(or {\it quantum stationary states})
in the case when
the Maxwell external potentials do not depend on time:
in this case, the Schr\"odinger equation 
reads as
\begin{equation}\label {Sc0}
i \hbar \dot \psi(x,t) =
H \psi (t): =
\frac{1}{2 \cm} 
\Big[- i \hbar \nabla- \displaystyle \frac e c \mathbf{A}_{\mathrm{ext}}(x)\Big]^2 \psi(x,t) + e A_{\mathrm{ext}}^0(x) \psi(x,t),
\end{equation}
where $m$ and $e$ denote the electron's mass and charge, $c$ is the speed of light
in vacuum,
$\mathbf{A}_{\mathrm{ext}}(x)$ is the external magnetic potential, and $A_{\mathrm{ext}}^0(x)$ is the external scalar potential.
In the {\it unrationalized}\,\footnote{\emph{Rationalization}
refers to
the choice of units aimed at
suppression of the explicit factors of $4\pi$
in the Maxwell equations.}
Gaussian units,
also called the Heaviside--Lorentz units
\cite[p. 781] {Jackson}, 
the values of the physical constants
(the electron charge\index{electron! charge} and mass\index{electron! mass},
the speed of light\index{speed of light} in  vacuum,
 the Boltzmann constant\index{Boltzmann! constant}, and the Planck constant\index{Planck! constant})
 are approximately equal to
\begin{equation} \label{HL}
\left.
\begin{array}{c}
e=- 4.8 \times 10^{-10} {\mathrm{esu}}, \qquad \cm=9.1 \times 10^{-28}
{\mathrm{g}},
\quad c=3.0 \times 10^{10} {\mathrm{cm/s}}
\\\\
k=1.38\times 10^{-16} {\mathrm{erg/K}},\qquad\,\,
\hbar=1.1 \times 10^{-27} \mathrm{erg}\cdot\mathrm{s}
\end{array}\right|.
\end{equation}

In Schr\"odinger's theory \cite{S1926},
stationary orbits
are identified
with finite energy solutions of the form
\begin{equation}\label{SSO}
\psi(x,t) = \varphi_\omega(x) e^{- i \omega t}, \qquad \omega \in \R.
\end{equation}
Substitution
of this Ansatz
into the Schr\"odinger equation (\ref{SSO}) leads to the famous eigenvalue problem
\begin{equation}\label{eip}
\omega\varphi_\omega=H\varphi_\omega,
\end{equation}
which determines the corresponding frequencies $\omega$ and amplitudes $\varphi_\omega$;
the latter are subject to the normalization condition
\begin{equation}\label{nor}
\int|\varphi_\omega(x)^2|\,dx=1.
\end{equation}
This condition means that the wave function describes one electron, as discussed below.
\smallskip

Such definition is rather natural, since then $ | \psi(x,t) | $ does not depend on time. Most likely, this definition was suggested
by the de Broglie wave function
 for {\it free particles} $ \psi(x,t) = Ce^{i (k x- \omega t)} $, which factorizes as $ Ce^{i k x} e^{- i \omega t} $. Indeed, in the case of  {\it bound particles},
it is natural to change the spatial factor $ Ce^{i k x} $, since the spatial properties have changed and ceased to be homogeneous.
On the other hand, the homogeneous time factor $ e^{- i \omega t} $ must be preserved, since the external potentials are independent of time.
However, these ``algebraic'' arguments do not 
withdraw the question of agreement of the Schr\"odinger definition
with the  Bohr postulate (\ref{B1})!
\smallskip

Thus, the problem of the mathematical interpretation of the Bohr postulate (\ref{B1}) in the Schr\"odinger theory  arises. 
   One of the simplest interpretations of the jumps (\ref{B1}),
   (\ref{B21})
is the long-time  asymptotics 
\begin{equation}\label{BS}
\psi(x,t) \sim \psi_\pm(x) e^{- i \omega_\pm t}, \qquad t \to \pm \infty
\end{equation}
{\it for each finite energy solution},
where
$ \omega_- = \omega_n $ and $ \omega_+ = \omega_{n'} $. However, 
{\it for  the linear Schr\"odinger equation} (\ref{Sc0}),
such asymptotics are obviously  wrong due to the {\it superposition principle}:
for example, such asymptotics
can not hold
for solutions of the form $ \psi(x,t) \equiv \phi_1(x) e^{- i \omega_1 t} + \phi_2(x) e^{- i \omega_2 t} $ with $\omega_1\ne\omega_2$
and with nonzero $\phi_1$ and $\phi_2$.
It is exactly this contradiction that shows that the linear 
Schr\"odinger equation alone cannot  serve as the basis for the 
theory compatible
with the Bohr postulates.

\begin{remark}
The Schr{\"o}dinger definition (\ref{SSO}) has been formalized by
P. Dirac \cite{DiracPQM} and
J. von Neumann
\index{von Neumann}
\cite{JvN}: they have identified  quantum states with the rays in the Hilbert space,
so the entire circular orbit (\ref{SSO}) corresponds to one quantum state.\index{quantum! state}
Our main conjecture is that this
formalization
is indeed a {\it theorem}: namely, we suggest that the
asymptotics\index{asymptotics} (\ref{BS}) is an inherent property of the nonlinear Maxwell--Schr{\"o}dinger system
(see \eqref{M}--\eqref{S} below). This
conjecture
is confirmed by perturbative arguments below.
\end{remark}

\section{Coupled Maxwell--Schr\"odinger equations}

As we have seen, the linear Schr\"odinger equation cannot explain the 
Bohr postulates. Thus, one should look for a nonlinear theory. 
Fortunately, 
we do not need to invent anything artificial since
such a theory is well established since 1926: it is the system
of the Schr\"odinger equation coupled to the Maxwell equations
which was essentially introduced in the first of
Schr\"odinger's articles \cite{S1926}. 
The corresponding second-quantized system 
is the main subject of nonrelativistic Quantum Electrodynamics.

This coupled system arises in the following way.
The wave function $\psi(x,t)$ defines the corresponding charge and 
current densities, which, in turn, generate their own Maxwell field.
The corresponding potentials
of this own field should be added to the external Maxwell
potentials in the Schr\"odinger equation (\ref{Sc0}).
This modification implies the nonlinear interaction between the wave function and the 
Maxwell field.
More precisely,
the Schr\"odinger theory associates to the wave function
$\psi(x,t)$ the corresponding charge and current densities
\begin{equation}\label{rj}
\rho(x,t) = e | \psi(x,t) |^2, \qquad \mathbf{j}(x,t)
=
\frac e \cm \Re
\Big(
\overline \psi(x,t) \Big[
- i \hbar \nabla- \frac ec \mathbf{A}_{\mathrm{ext}} (x,t)
\Big] \psi(x,t)
\Big).
\end{equation}
For any solution of the Schr\"odinger equation (\ref{Sc0}),
these densities satisfy the charge continuity equation
\begin{equation}\label{cce}
\dot\rho(x,t)+\dv \mathbf{j}(x,t)=0.
\end{equation}
Note that the normalization condition 
(\ref{nor}) means that
the total charge corresponding to the stationary orbit (\ref{SSO}) is
\begin{equation}\label{nor2}
\int\rho(x)\,dx=e,
\end{equation}
i.e., we have exactly one electron on this orbit.
The charge and current densities (\ref{rj}) generate the Maxwell field according
to the Maxwell equations
\begin{equation}\label{mhl-1}
\begin{cases}
\dv \mathbf{E}(x,t)=\rho(x,t),
\quad
&\curl \mathbf{E}(x,t)=- \frac 1c \dot{\mathbf{B}}(x,t),
\\[2ex]
\dv \mathbf{B}(x,t)=0,
\quad
&\curl \mathbf{B}(x,t)=\frac 1c
\big(\mathbf{j}(x,t) + \dot{\mathbf{E}}(x,t)\big).
\end{cases}
\end{equation}
The second and third equations imply the Maxwell representations
\begin{equation}\label{EBA}
 \mathbf{B}(x,t)=\curl \mathbf{A}(x,t),\qquad
 \mathbf{E}(x,t)=- \frac 1c \dot{\mathbf{A}}(x,t) -\nabla A^0(x,t).
 \end{equation}
We will assume the Coulomb gauge
\begin{equation}\label{cgauge2}
  \dv  \mathbf{A}(x,t) \equiv 0.
\end{equation}
Then  
the Maxwell equations\index{Maxwell! equations} (\ref{mhl-1}) are equivalent to the system
\begin{equation}\label{M}
\frac 1 {c^2} \ddot{\mathbf{A}}(x,t)=\Delta \mathbf{A}(x,t) +\frac 1c P \mathbf{j}(x,t),
\quad \Delta A^0(x,t)=-\rho(x,t),\qquad x \in\R^3,
 \end{equation}
where $P$ is the \emph{orthogonal projection}
in the real Hilbert space $L^2(\R^3)\otimes\R^3$
onto divergence-free vector fields.
\smallskip

Finally, one should add these ``own'' Maxwell potentials $\mathbf{A}(x,t)$ and $A^0(x,t)$ 
to the external potentials in the Schr\"odinger equation (\ref{Sc0}), so we obtain
the modified Schr\"odinger equation
\begin{eqnarray}\label {S}
i \hbar \dot \psi(x,t) &=&
H(t) \psi (t)
\nonumber\\
\nonumber\\
&: =&
\frac{1}{2 \cm} 
\Big[- i \hbar \nabla- \displaystyle \frac e c (\mathbf{A}_{\mathrm{ext}}(x)+\mathbf{A}(x,t))\Big]^2 \psi(x,t) 
+ e (A_{\mathrm{ext}}^0(x)+A^0(x,t)) \psi(x,t).
\qquad\qquad
\end{eqnarray}
Equations (\ref{M})--(\ref{S}) form the nonlinear
Maxwell--Schr\"odinger system.

\section{Bohr's postulates via perturbation theory}\label{sBpt}   
   
The remarkable success of the Schr\"odinger theory was in the explanation
of  Bohr's postulates via asymptotics (\ref{BS})
 by 
 means of 
{\it perturbation theory} applied to the 
{\it coupled Maxwell--Schr\"odinger equations} (\ref{M})--(\ref{S})
in the case of {\it static external potentials}
\begin{equation}\label{stat}
 \mathbf{A}_{\mathrm{ext}}(x,t)\equiv \mathbf{A}_{\mathrm{ext}}(x),\qquad
 A^0_{\mathrm{ext}}(x,t)\equiv A^0_{\mathrm{ext}}(x).
 \end{equation}
 For instance, in the case of an atom,  $A^0_{\mathrm{ext}}(x)$ is the Coulomb potential of the nucleus, while
 $\mathbf{A}_{\mathrm{ext}}(x)$ is  the vector  potential of the magnetic field of the nucleus.
Namely, as the first approximation, the fields $\mathbf{A}(x,t) $ and $ A^0(x,t) $ in equation 
(\ref{S}) can be neglected, so we obtain 
equation (\ref{Sc0}) with the Schr\"odinger operator $H(t)\equiv H$.
For 
``sufficiently good''
external potentials and initial conditions,
each finite energy solution
can be expanded into eigenfunctions:
\begin{equation}\label{Sexp}
\psi(x,t) = \sum_n C_n \psi_n(x) e^{- i \omega_n t} + \psi_c(x,t), \qquad
\psi_c(x,t) =
\int C (\omega) e^{- i \omega t}\,d \omega,
\end{equation}
where the integration is performed over the continuous spectrum of the operator $ H $, and 
 for any $R>0$ we have
\begin{equation}\label{psic}
\int_{|x|<R}|\psi_c(x,t)|^2\,dx\to 0,\quad t\to\pm\infty;
\end{equation}
see, for example, \cite[Theorem 21.1] {KopK2012}.
The
substitution of the expansion (\ref{Sexp}) into the expression for currents (\ref{rj}) gives 
\begin{equation}\label{jexp}
\mathbf{j}(x,t) = \sum_{n,\,n'} \mathbf{j}_{nn'}(x) e^{- i \omega_{nn'} t} + c.\,c. + \mathbf{j}_c(x,t),
\end{equation}
where $ \mathbf{j}_c(x,t) $ has a continuous frequency spectrum. 
Thus, the current in the right-hand side of 
the
Maxwell equation  (\ref{M}) contains,
besides the continuous spectrum, only
{\it discrete frequencies} $ \omega_{nn'} $.
Hence, the {\it discrete spectrum} of the corresponding 
Maxwell field also contains only these frequencies $ \omega_{nn'} $.
This justifies the Bohr rule (\ref{B21})
{\it only in the first order of perturbation theory},
since this calculation
ignores the feedback action of the radiation onto the atom.
\smallskip

The same arguments also suggest treating
 the jumps (\ref{B1}) as
the {\it single-frequency asymptotics} (\ref{BS})
for solutions of the  
 Maxwell--Schr\"odinger system (\ref{M})--(\ref{S}).
Indeed, the currents (\ref{jexp}) on the right of
the Maxwell equation (\ref{M}) produce radiation when
nonzero frequencies $ \omega_{nn'}\ne 0 $ are present. 
This is due to the fact that $\R\setminus 0$ is the
absolutely  continuous spectrum of the Maxwell
equations.
However, this radiation
cannot last forever, since
it
irrevocably
carries the energy to infinity while
the total energy is finite. Therefore, in the long-time
limit only $ \omega_{nn'}=0 $ survives, which means exactly 
that we have
 {\it single-frequency asymptotics} (\ref{BS})
in view of  (\ref{psic}).

\section{Quantum jumps as global attraction}
Of course, the perturbation arguments above cannot provide a rigorous justification
of the long-time asymptotics (\ref{BS}) for
the coupled Maxwell--Schr\"odinger equations. 
We suggest that this asymptotics 
 is a manifestation of a more general fact which we will discuss in the next section.
 Note that the system (\ref{M})--(\ref{S}) admits the symmetry group
 $\mathbf{U}(1)$:
\begin{equation}\label{U1-2}
(\mathbf{A}(x),   A^0(x), \psi(x)) \mapsto (\mathbf{A}(x),  A^0(x), \psi(x) e^{ i \theta}), \qquad \theta \in [0,2 \pi].
\end{equation}
That is, if $(\mathbf{A}(x,t), A^0(x,t), \psi(x,t))$
is a solution
to the system (\ref{M})--(\ref{S}),
then,
for any $\theta\in\R$,
$(\mathbf{A}(x,t), A^0(x,t), \psi(x,t) e^{ i \theta})$
is also a solution.

\begin{definition}
{\it Stationary orbits} of the Maxwell--Schr\"odinger nonlinear system
(\ref{M})--(\ref{S}) are
finite energy
solutions of  the form 
\begin{equation}\label{Sorb}
\big(\mathbf{A}(x),\,A^0(x),\, e^{- i \omega t} \phi(x)\big).
\end{equation}
\end{definition}

The existence of such stationary orbits for the system (\ref{M}), (\ref{S})
 was proved in \cite{CG2004} in the case of external potentials
\begin{equation}\label{AAA}
\mathbf{A}_{\mathrm{ext}}(x,t) \equiv 0, \qquad A_{\mathrm{ext}}^0(x,t) =-\frac {e Z} {| x |}.
\end{equation}
The perturbation arguments of previous section suggest that the Bohr postulate (\ref{B1}) can be treated as
 the long-time asymptotics 
\begin{equation}\label{ate1}
\big(\mathbf{A}(x,t),\,A^0(x,t),\,\psi(x,t)\big)
\sim \big(\mathbf{A}_\pm(x),\,A^0(x),\,e^{- i \omega_\pm t} \phi_\pm(x)\big),
\qquad t \to \pm \infty
\end{equation}
{\it for all finite-energy solutions} of the  Maxwell--Schr\"odinger equations
(\ref{M}),  (\ref{S}) in the case of static external potentials (\ref{stat}). We conjecture that these asymptotics 
hold 
in the $H^1$-norm on every bounded region  of $\R^3$.
The asymptotics (\ref{ate1}) 
 means that there is a {\it global attraction} to
the set of stationary orbits.
We expect
that a 
similar attraction takes place for
Maxwell--Dirac, Maxwell--Yang--Mills, and other
coupled equations.
This suggests the interpretation of quantum stationary
states  as the points and orbits that constitute the {\it global attractor} of the corresponding quantum dynamical equations.

The asymptotics  (\ref{ate1}) have not been proved
 for the Maxwell--Schr\"odinger system
(\ref{M})--(\ref{S}). On the other hand, similar asymptotics
are now being proved 
in \cite{C2012}--\cite{KK2019} for a number of 
 model nonlinear PDEs with the symmetry group $\mathbf{U}(1)$.
 In the 
next section we state a general conjecture which 
reduces to the asymptotics (\ref{ate1}) in the case of the Maxwell--Schr\"odinger system.

\begin{remark}
Experiments show that the transition time of quantum jumps (\ref{B1}) is of the order of $10^{-8}$s, although the asymptotics  (\ref{ate1}) require infinite time.
We suppose that this discrepancy can be explained by the following arguments:
\smallskip

\noindent
i) $10^{-8}$s is the transition  time between very small neighborhoods of
initial and final states;
\smallskip

\noindent
ii)
during this time,
the atom emits the major part of the radiated energy.
\end{remark}

\section{Relation to linear Quantum Mechanics}
The  linear eigenvalue problem (\ref{eip})
together with
Bohr's second rule (\ref{B21}) give very efficient
description of the atomic and molecular spectra.
This efficiency
is presumably
due to the smallness
of the fraction $\frac ec\sim 10^{-20}$.
Namely,
substituting the  stationary orbit  (\ref{Sorb}) into the 
 Maxwell--Schr\"odinger system (\ref{M})--(\ref{S}), 
 we obtain the corresponding \emph{nonlinear eigenvalue problem}
for the  frequencies $\omega$, which
play the role of \emph{nonlinear eigenvalues},
and for the corresponding \emph{nonlinear eigenfunctions}
$(\mathbf{A}(x),A^0(x),\phi(x))$.
After such a substitution,
we can neglect the terms with
$\mathbf{A}(x)$ and $A^0(x)$ in equation (\ref{S});
this corresponds to the first order approximation
of the perturbation series in $\frac ec$.
Then the nonlinear eigenvalue problem
reduces to the linear one, represented by (\ref{eip}).

\section{Conjecture on attractors of $G$-invariant equations}

In this section, we formulate the general conjecture on attractors 
introduced in  
\cite{K2016,KK2020,KK2021}.
Let us consider general 
 $G$-invariant  {\it autonomous} Hamiltonian nonlinear PDEs in $\R^n$ 
 of the form
 \begin{equation}\label{dyn}
\dot\Psi(t)=F(\Psi(t)),\qquad t\in\R,
\end{equation}
with a Lie  symmetry group $G$ acting
on a suitable Hilbert or Banach phase space $\mathcal{E}$ via a linear 
representation $T$.
The Hamiltonian structure  means that
\begin{equation}\label{Hamstr}
F(\Psi)=J D \mathcal{H}(\Psi),\qquad J^*=-J,
\end{equation}
where $\mathcal{H}$ denotes the corresponding Hamiltonian functional.
The
$G$-invariance means that 
\begin{equation}\label{Ginv}
 F (T(g) \Psi) = T(g)F (\Psi),\qquad
 \Psi\in \mathcal{E},
\end{equation}
for all $g \in G$.
In that case, for any solution $ \Psi (t) $ to
equation (\ref{dyn}), the trajectory  $ T(g) \Psi (t) $ is also a solution,
so the representation commutes with the dynamical group 
$U(t):\Psi(0)\mapsto \Psi(t)$; that is,
\begin{equation}\label{comTU}
T(g)U(t)=U(t)T(g).
\end{equation}
Let us note that
the theories of elementary particles deal systematically with the
symmetry groups
$\mathbf{SU}(2)$, $\mathbf{SU}(3)$, $\mathbf{SU}(5)$, $\mathbf{SO}(10)$
and so on,
like e.g. the group
$$
\mathbf{SU}(4) \times \mathbf{SU}(2) \times \mathbf{SU}(2)
$$
from the Pati--Salam model \cite{AF1963,PS1974},
one of the candidates
for the ``Grand Unified Theory''.

\smallskip

\noindent
{\bf Conjecture  A (on attractors).} {\it
For ``generic'' $G$-invariant autonomous equations (\ref{dyn}),
any finite energy solution $ \Psi (t) $  admits a long-time asymptotics
\begin {equation} \label {at10}
\Psi (t) \sim e^{\hat \lambda_\pm t} \Psi_\pm, \qquad t \to \pm \infty,
\end {equation}
in the appropriate topology of the phase space $\mathcal{E}$.
Here 
$\hat \lambda_{\pm}=T'(e)\lambda_{\pm}$, where
$e\in G$ is the identity element and
$\lambda_{\pm}$ belong to the corresponding 
 Lie algebra
$\mathfrak{g}$, 
while $\Psi_{\pm}$ are some
\emph{limiting amplitudes} depending on the solution $\Psi(t)$.
}
\smallskip

In other words, all solutions of the form
$e^{\hat\lambda t} \Psi$ with $\hat\lambda=T'(e)\lambda$,
where
$\lambda\in\mathfrak{g}$,
constitute
a global attractor for \emph{generic} $G$-invariant Hamiltonian nonlinear
PDEs of the form \eqref{dyn}.
This conjecture suggests
that we define {\it stationary  $G$-orbits} for equation \eqref{dyn} as solutions of the form
\begin {equation} \label {at10a}
\Psi (t) = e^{\hat \lambda t} \Psi,
\qquad t\in\R;
\qquad
\hat\lambda=T'(e)\lambda,\quad\lambda\in\mathfrak{g}.
\end {equation}
 This definition leads to the corresponding
 {\it nonlinear eigenvalue problem}
\begin {equation} \label {at10b}
\hat \lambda \Psi = F (\Psi).
\end {equation}
In particular,
in the case of the linear Schr\"odinger equation
with the symmetry  group  $\mathbf{U}(1) $,
stationary orbits are solutions of the form $ e^{i \omega t} \phi(x) $,
where $ \omega \in \R $ is an eigenvalue of the Schr\"odinger operator
and $ \phi(x) $ is the corresponding eigenfunction. 
In the case of the symmetry group $G=\mathbf{SU}(3)$, the generator
(``eigenvalue'') $ \lambda $ is
a $3\times 3$-matrix, and solutions (\ref{at10a})
can be, in particular, quasiperiodic in time.
\smallskip

Note that
the conjecture  (\ref{at10}) fails for linear equations:
that is, linear equations are exceptional, not ``generic''!

\begin{remark}
The existence of the stationary orbits of the form (\ref{at10a})
was proved  for nonlinear wave equations and the Maxwell--Schr\"odinger
equations \cite{BL83-1,BL83-2,CG2004,St77}. The proofs of the existence
for the Maxwell--Dirac and Klein--Gordon--Dirac
coupled nonlinear equations in \cite{EGS} essentially 
relies on the topological methods of Lusternik--Schnirelman \cite{LS1934,LS1947}.
\end{remark}

\medskip

\noindent
{\bf Empirical evidence.}
The conjecture (\ref{at10}) agrees with the Gell-Mann--Ne'eman theory of baryons  \cite{GM1962,Ne1962}.
Indeed, 
in 1961,
Gell-Mann and Ne'eman suggested
using the symmetry group $\mathbf{SU}(3)$ 
for the strong interaction of baryons
relying on the discovered parallelism between empirical data
for the baryons
and the ``Dynkin diagram''
of the Lie algebra
$\mathfrak{g}=\mathfrak{su}(3)$ with $8$ generators (the famous ``eightfold way'').
This theory resulted in
the scheme of quarks in
quantum chromodynamics \cite{HM1984}
and in the prediction of a new baryon
with prescribed values of
its mass and decay products. This particle (the $\Omega^-$-hyperon)
was promptly discovered experimentally \cite{omega-1964}.

On the other hand, the elementary particles seem to describe long-time asymptotics
of quantum fields. Hence
the empirical correspondence between elementary
particles and generators of the Lie algebras
presumably gives an
evidence in favor of our general conjecture
(\ref{at10}) for equations with  Lie symmetry groups.

\section {Results on global attractors for nonlinear Hamiltonian PDEs}
\label{sattr}

Here we give a  brief survey of
rigorous results \cite{K1995a}--\cite{C2013}
obtained since 1990 that confirm the conjecture  (\ref{at10})
for a list of model equations of the form (\ref{dyn}).
The details can be found in the surveys \cite{K2016,KK2020,KK2021}.

 The results 
confirm the existence of finite-dimensional attractors
in the Hilbert or Banach phase spaces and demonstrate
the explicit correspondence between
the long-time asymptotics and  the symmetry group $G$ of equations.
\smallskip
The results obtained so far 
concern equations (\ref{dyn}) with
the following
four basic groups of symmetry:
the trivial symmetry group
$G=\{e\}$, the translation group $G=\R^n$ for 
translation-invariant equations, the unitary group $G=\mathbf{U}(1)$
for phase-invariant equations,
and the orthogonal group $\mathbf{SO}(3)$ for ``isotropic'' equations.
In these cases, the asymptotics (\ref{at10}) reads as follows.

\subsection{Equations with the trivial symmetry group $G=\{e\}$}
\begin{figure}[htbp]
		\begin{center}
			\includegraphics[width=0.7\columnwidth]{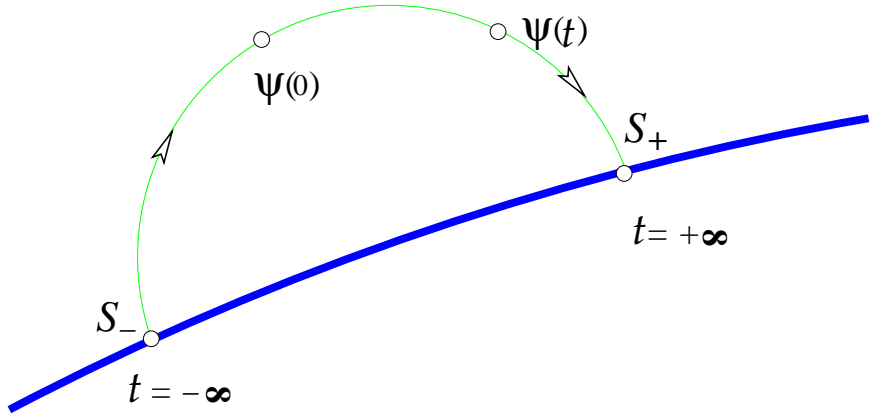}
			\caption{Convergence to stationary states.}
			\label{fig-1}
		\end{center}
	\end{figure}

For such {\it generic equations},
the conjecture (\ref{at10})  means
the  attraction  of
 {\it all finite energy solutions}
 to  stationary states:
 \begin {equation} \label {ate}
\psi(x,t) \to S_\pm(x), \qquad t \to \pm \infty,
\end {equation}
as illustrated on Fig.~\ref{fig-1}.
Here  the states $S_{\pm}=S_{\pm}(x) $  depend on  the trajectory $\Psi(t)$ under
consideration,
while the convergence holds in local 
seminorms of type  $L^2(|x|<R)$ with any $R>0$.
This convergence cannot hold
in global norms (i.e., in norms corresponding to $R=\infty$)  due to  
energy conservation.
The asymptotics (\ref{ate}) can be symbolically written as the transitions
\begin{equation}\label{Bohr}
S_-\mapsto S_+,
\end{equation}
which
can be considered as the mathematical model of
Bohr's ``quantum jumps'' (\ref{B1}).
Such an attraction
was proved for a variety of model equations
in \cite{D2012}--\cite{T2010}.
In \cite{K1991}--\cite{KSK1997}, the convergence was proved
i) for a string coupled to nonlinear oscillators:
\begin{equation}\label{sos}
\ddot\psi(x,t)=\psi''(x,t)+f(x,\psi(x,t)),\qquad x\in\R,
\end{equation}
with $f(x,\psi(x,t))=\sum_k\delta(x-x_k)F_k(\psi(x_k,t))$
and with $f(x,\psi(x,t))=\chi(x)F(\psi(x,t))$;
ii) for a~three-dimensional wave equation coupled
to a charged particle 
\begin{equation}\label{wp}
\left\{\begin{array}{rcl}
\ddot\psi(x,t)=\Delta\psi(x,t)-\rho(x-q(t)),\qquad x\in\R^3
\\[1ex]
\dot p(t)=-\nabla V(q)-\displaystyle\int\nabla\psi(x,t)\rho(x-q(t))\,dx,
\end{array}
\right|
\end{equation}
where
\[
p(t)=\frac{m\dot q(t)}{\sqrt{1-\dot q^2(t)}}
\]
is the relativistic momentum
of the particle. This is Hamiltonian system with the Hamilton functional
\begin{equation}
\label{Ham}
H(\psi,\pi,q, p)=\frac12\int[|\pi(x)|^2+|\nabla\psi(x)|^2]\,dx
+\int\psi(x)\rho(x-q)\,dx
+\sqrt{1+p^2}+V(q)
\end{equation}
The global attraction to solitons  was extended i) to 
the wave equation, the Dirac equation,
and the Klein--Gordon equation
with concentrated nonlinearities, and
 ii) to the Maxwell--Lorentz system,
which is the system of type (\ref{wp})
with the Maxwell equations instead of the wave equation
and the Lorentz equation
instead of the Newton equation:
\begin{equation} \label{mls3}
\left\{
\begin{array}{l}
\dot{\mathbf{E}}(x,t)=\curl \mathbf{B}(x,t)-\dot q(t)\rho(x-q(t)),
\qquad\qquad
\dot{\mathbf{B}}(x,t)= - \curl \mathbf{E}(x,t)
\\[1.5ex]
\dv \mathbf{E}(x,t)=\rho (x-q(t)),
\qquad\qquad\qquad\quad \quad \quad\quad
\dv \mathbf{B}(x,t)=0
\\[1.5ex]
\dot p(t)=-\nabla V(q)+\displaystyle\int [\mathbf{E}(x,t)+\dot q(t)\times\mathbf{B}(x,t)]\rho(x-q(t))\,dx
\end{array}
\right|.
\end{equation}
This model of the electrodynamics 
with {\it extended electron}
was introduced by Abraham to avoid the infinite mass and energy
of the charged point particle,
known as the {\it ultraviolet divergence}, \cite{A1902,A1905}.

All proofs rely on the analysis of radiation 
which  irreversibly carries a portion of  energy
to infinity.
The  details can be found  in the surveys \cite{KK2020,KK2021}. 

In all the problems considered,
the convergence (\ref{ate}) implies, by the Fatou theorem, the 
inequality
\begin{equation}\label{1.6''}
{\cal H}(S_{\pm})\leq {\cal H}(Y(t))\equiv\const,\,\,\,t\in\R,
\end{equation}
where $\mathcal{H}$ is the corresponding Hamiltonian (energy) 
functional. This inequality is an analog
 of the well-known property of weak convergence in 
Hilbert and Banach spaces. 
Simple examples show that
the strict inequality in (\ref{1.6''}) is possible, which
means that an irreversible scattering of energy  to infinity occurs.

\begin{example}[The d'Alembert waves]
{\rm
In particular, the asymptotics (\ref{ate})
with the strict inequality (\ref{1.6''}) 
can 
easily be demonstrated for the d'Alembert equation 
$\ddot\psi(x,t)=\partial_x^2\psi(x,t)$
with general solution
\begin{equation}\label{dalw}
\psi(x,t)=f(x-t)+g(x+t).
\end{equation}
 Indeed, the  convergence 
$\psi(\cdot,t)\to 0$ in $L^2_{\mathrm{loc}}(\R)$ obviously
holds for all $f,g\in L^2(\R)$. On the other hand, 
 the convergence to zero in {\it global norms}
obviously fails if $f(x)\not\equiv 0$ or $g(x)\not\equiv 0$. 
}
\end{example}

\begin{example}[Nonlinear Huygens Principle]
{\rm
Consider solutions of  3D wave equation  with a unit propagation velocity and initial data with support in a ball $ | x | <R $. The corresponding
solution is concentrated in spherical layers $ | t | -R <| x | <| t | + R $.
Therefore,
the energy localized 
in any bounded region converges  to zero
as
$ t \to \pm \infty $, although its total energy remains constant.
This convergence to zero is  known as the
{\it  strong  Huygens principle}. Thus, global attraction to stationary states
(\ref{ate}) is a generalization of the strong Huygens principle to nonlinear equations.
The difference is that for the linear wave equation, the limit  
 is always zero, while for nonlinear equations
the limit can be any stationary solution.
}
\end{example}

\begin{remark}
The proofs
in \cite{KSK1997} and \cite{KS2000}
rely on the relaxation of acceleration
\begin{equation}\label{rel}
\ddot q (t) \to 0, \qquad t \to \pm \infty
\end{equation}
for solutions to the coupled  system (\ref{wp})
and similarly for solutions to the Maxwell--Lorentz system.
This relaxation was discovered about 100 years ago
in classical electrodynamics
and is known as the \emph{radiation damping};
see \cite[Chapter 16]{Jackson}.
The detailed account on early investigations of this problem
 by A. Sommerfeld, H.~Poincar\'e, P. Dirac and others
 can be found in Chapter 3 of \cite{S2004}, see also \cite{Dirac1938}.

However, a~rigorous proof of the relaxation\index{relaxation of acceleration} (\ref{rel})
is not so obvious.  
It first appeared in \cite{KSK1997} and \cite{KS2000} under the
Wiener condition  on the particle charge density
\begin {equation} \label {W1}
\hat \rho (k): = \int e^{ikx} \rho(x)\,dx \ne 0, \qquad k \in \R^3.
\end {equation}
This condition is an
analogue of the ``Fermi Golden Rule''
first introduced by Sigal in the context of nonlinear wave and Schr\"odinger equations \cite{Sig1993}.
The proof of the relaxation (\ref{rel}) relies
on a novel application of the Wiener Tauberian theorem.
\end{remark}

\subsection{Group of translations $G=\R^n$}

For  {\it generic translation-invariant equations},
the conjecture (\ref{at10}) means
the attraction of {\it all finite-energy solutions} to solitons:
\begin {equation} \label {att}
\psi(x,t) \sim \psi_\pm (x-v_\pm t), \qquad t \to \pm \infty,
\end {equation}
where the
convergence holds in local seminorms
of type $ L^2 (| x-v_\pm t | <R) $ with any $ R> 0 $, i.e.,
{\it in the comoving reference frame}.
A trivial example is provided by  the d'Alembert equation
$ \ddot \psi(x,t) = \partial_x^2\psi(x,t) $ with 
general solution 
 (\ref{dalw}) corresponding to the asymptotics 
(\ref{att}) with $v_+=\pm 1$ and $v_-=\pm 1$.
\smallskip 

Such soliton asymptotics 
was first proved   for {\it completely integrable equations}
(Korteweg--de Vries equation, etc.); see \cite{EvH,Lamb80}.
Moreover, for  the Korteweg--de Vries equation,
more accurate soliton asymptotics in {\it global norms}
with several solitons were discovered
by Kruskal and Zabusky in 1965 by
numerical simulation:
it is the decay to solitons
\begin{equation}\label{attN}
\psi(x,t) \sim \sum_k \psi_\pm (x-v^k_\pm t) + w_\pm(x,t), \qquad t \to \pm \infty,
\end{equation}
where
$ w_\pm $ are some dispersive waves.
\smallskip 

 Later on, such asymptotics
were proved by the method of {\it inverse scattering problem}
for nonlinear
completely integrable Hamiltonian
translation-invariant equations (Korteweg--de Vries, etc.)
in the works of Ablowitz, Segur, Eckhaus, van Harten, and others \cite{EvH,Lamb80}.
\smallskip

For equations which are not completely integrable,
the global attraction (\ref{att}) was established 
in \cite{IKM2004}--\cite{KMV2004}.
The first result was obtained in
\cite{KS1998}
for the translation-invariant system (\ref{wp}) with zero external potential $V=0$.
In this case the total momentum is conserved:
\begin{equation}\label{PP}
P:=p-\int \dot\psi(x,t)\nabla\psi(x,t)\,dx=\const.
\end{equation}
In
\cite{IKM2004}, the result was extended to
the Maxwell--Lorentz system (\ref{mls3}) with $V=0$.
The proofs
in \cite{KS1998} and \cite{IKM2004}
rely on  variational properties of solitons and their orbital stability, as well as on the  relaxation of the acceleration (\ref{rel})
under the Wiener condition  (\ref{W1}).
The case of small charge was considered in \cite{IKM2003, IKS2004a}.

Let us mention that
for non-completely-integrable equations
the multi-soliton asymptotics 
(\ref{attN}) 
were observed numerically 
for 1D {\it relativistically-invariant} nonlinear wave equations
in \cite{KMV2004}.

\subsection{Unitary symmetry group $G=\mathbf{U}(1)$}

\begin{figure}[htbp]
	\begin{center}
		\includegraphics[width=0.9\columnwidth]{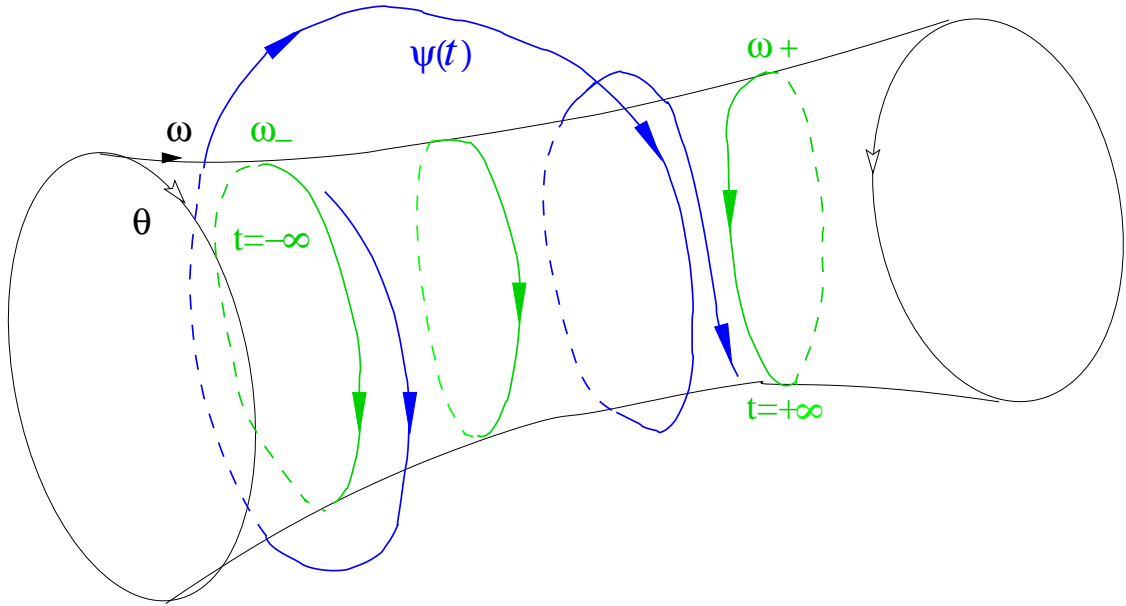}
		\caption{Convergence to stationary orbits.}
		\label{fig-3}
	\end{center}
\end{figure}

For {\it generic $\mathbf{U}(1)$-invariant equations},
the conjecture (\ref{at10}) means
the  attraction of {\it all finite-energy solutions}
to ``stationary orbits'':
\begin{equation}\label{atU}
\psi(x,t)\sim\phi_\pm(x) e^{-i \omega_{\pm} t}, \qquad t \to \pm \infty,
\end{equation}
where
 $\omega_{\pm}\in\R$.
 Such  asymptotics are similar to 
  the Bohr transitions  between stationary orbits  (\ref{ate1}) of the coupled Maxwell--Schr\"odinger equations. 
This asymptotics means that there is a
global attraction to the solitary manifold formed by all
\emph{stationary orbits} (\ref{SSO}).
The asymptotics is considered
in the local seminorms $L^2(|x|<R)$ with any $R>0$.
The global attractor is a smooth manifold formed by the circles
that are the orbits of the action of the 
symmetry group $\mathbf{U}(1)$, as is illustrated on
Fig.~\ref{fig-3}.

Such an attraction 
  was proved 
  for the first time
  i) in \cite{C2012} and
\cite{K2003}--\cite{KK2010b}
for the 
Klein--Gordon
  and Dirac equations coupled to
$\mathbf{U}(1)$-invariant nonlinear oscillators,
\begin {eqnarray}
\ddot \psi(x,t) & = & \Delta \psi(x,t) -m^2 \psi +
\sum_{k = 1}^N \rho (x - x_k) F_k (\langle \psi (\cdot, t), \rho (\cdot -x_k) \rangle), \label {KGn} \\
\nonumber \\
i \dot \psi(x,t)
& = &
\big (-i \alpha \cdot \nabla + \beta m \big) \psi
+ 
\sum_{k = 1}^N \rho (x - x_k) F_k (\langle \psi (\cdot, t), \rho (\cdot -x_k) \rangle),
 \label{Dn}
\end {eqnarray}
under the Wiener condition\index{Wiener condition} \eqref {W1},
where $ \alpha = (\alpha_1, \dotsc, \alpha_n) $ and $ \beta = \alpha_0 $
are the Dirac matrices;
 ii) in
\cite{C2013} for a discrete  approximations
of such coupled systems,
i.e. for the corresponding finite difference schemes;
and iii) in \cite{K2017}--\cite{KK2019}
for the wave, Klein--Gordon, and Dirac equations with concentrated nonlinearities.
More precisely, we have proved  
global attraction to the 
{\it solitary manifold} of all 
stationary orbits, although
global attraction to particular  
stationary orbits (orbits with fixed $\omega_{\pm}$)
is still an open problem.

All
these results were proved under the assumption that the equations
are ``strictly nonlinear''. For linear equations, the global attraction obviously fails
if the discrete spectrum consists of at least two different eigenvalues.

The proofs of all these results  rely on i) a nonlinear analog of the Kato theorem on the absence of embedded eigenvalues;
ii) a new theory of multiplicators
in the space of quasimeasures;
iii) a novel application of the
Titchmarsh convolution theorem \cite{Hor90,Lev96,Tit26}.
For application to
the Klein--Gordon equation in discrete space-time (a numerical scheme)
in \cite{C2013},
the appropriate version of the Titchmarsh theorem is 
obtained in \cite{KK2013t}.

The asymptotics (\ref{atU})
is closely related to the {\it limiting amplitude principle}
\cite{Lad1957,Mor1962} and to the theory of waveguides \cite{Lewin}.

\subsection{Orthogonal group $G=\mathbf{SO}(3)$} 

For {\it generic rotation-invariant equations},
the conjecture (\ref{at10}) takes the form of
the long-time asymptotics
\begin{equation}\label{atSO3}
  \Psi(t)\sim  e^{-i \hat\Omega_{\pm} t}\Psi_{\pm}, \qquad t \to \pm \infty
\end{equation}
for {\it all finite-energy solutions}, 
where $\hat\Omega_{\pm}$ are 
suitable representations
of
$\Omega_{\pm}\in \mathfrak{so}(3)$.
This  means that
global attraction to ``stationary $\mathbf{SO}(3)$-orbits'' takes place.  
Such asymptotics
for the Maxwell--Lorentz equations with a rotating particle
are proved in \cite{IKS2004b}.

\section{\ On generic equations}
We must still    specify the 
meaning of the 
term \emph{generic} in  our conjecture (\ref{at10}).
 In fact, this conjecture
means that the asymptotics (\ref{at10}) hold for all solutions from
an open dense set of $G$-invariant equations.
\smallskip\\
i) In particular, the asymptotics
 (\ref{ate}), (\ref{att}), (\ref{atU}), and (\ref{atSO3})
hold under appropriate conditions, which define some ``open dense subset'' of
$G$-invariant equations with the four types of the symmetry group $G$.
This asymptotic expression
may break down if these conditions fail: this corresponds to some ``exceptional''
equations. For example, the global attraction  (\ref{atU}) breaks down for linear Schr\"odinger equations with at least two different eigenvalues.
Thus, linear equations are exceptional,
not  generic!
\smallskip\\
ii) The general situation is as follows. Let the Lie group $G_1$ be a (proper) subgroup of some larger Lie group $G_2$. Then
$G_2$-invariant equations
form an ``exceptional subset'' among all 
$G_1$-invariant equations, and
the corresponding asymptotics (\ref{at10}) may be completely different.
For example, the trivial group $ \{e \} $ is a proper 
subgroup in $\mathbf{U}(1) $ and in $ \R^n $, and the asymptotic expressions
(\ref{att}) and (\ref{atU}) may differ significantly from (\ref{ate}).
\smallskip\\
iii)
Examples
of solitary wave solutions
violating asymptotics (\ref{at10})
in particular \emph{non-generic} models
were constructed in
\cite{KK2009},
\cite{KK2010a}
(two-frequency solitary waves
in models based on the Klein--Gordon equation
with $\mathbf{U}(1)$-symmetry group),
in \cite{C2013}
(two- and four-frequency solitary waves
in the discrete time-space
Klein--Gordon equation with $\mathbf{U}(1)$-symmetry group),
and in \cite{BC2018}
(two-frequency solitary waves
in the Soler model with $\mathbf{SU}(1,1)$-symmetry group,
when the solitary manifold is larger than the orbit
of the action of $\mathbf{SU}(1,1)$).


\section{\ Adiabatic effective dynamics of solitons}

 The system \eqref {wp}
 admits soliton solutions 
 \begin{equation}\label{solsol}
 \phi_v(x-vt),\qquad q=vt
 \end{equation}
  in the case of identically
 zero external potential\index{external! potential} $V(q)\equiv 0$.
However, even in the case of
 nonzero external potential, \emph{soliton-like solutions} of the form
\begin {equation} \label{asol}
\psi(x,t) \approx \psi_{v (t)} (x-q (t))
\end {equation}
may exist if the potential is slowly varying:
\begin {equation} \label{asolV}
| \nabla V (q) | \le \varepsilon \ll 1.
\end {equation}
In this case, the  total momentum\index{conservation! of momentum} (\ref{PP}) is generally not conserved, but its slow evolution together with evolution
of the parameters $q(t)$, $v(t)$
in (\ref{asol}) can be described in terms of an appropriate finite-dimensional Hamiltonian dynamics.\index{Hamiltonian! equation}

Indeed, denote by $ P=P_v\in\R^3 $ the total momentum (\ref{PP}) of the soliton
 (\ref{solsol}).
It is important that the map $ \mathbf{P}: v \mapsto P_v $ is an isomorphism of the ball $ | v |< 1 $ on $ \R^3 $.
Therefore, we can consider $ q,\,P $ as global coordinates on the
\emph{solitary manifold} 
 $$
 \mathcal{S}:=\{S_{q,v}=(\phi_v(x-q),-v\cdot\nabla\phi_v(x-q),a,v):q,v\in\R^3,\,\,\,\,|v|<1\}.
 $$
The effective Hamiltonian functional is defined by
\begin {equation} \label{Heff}
H_{\rm eff} (q,P_v) \equiv H (S_{q, v}), \qquad q,\,P_v \in\R^3,
\end {equation}
where $H $ is the Hamiltonian (\ref{Ham})
with $V(q)\equiv 0$.
This functional
can be represented as
$$
H_{\rm eff} (Q, \Pi)=E (\Pi) + V (Q),
$$
since the
first integral in (\ref{Ham}) does not depend on $Q$, while the
last integral vanishes on the solitons.
 Hence, the corresponding Hamiltonian equations  have the form
\begin {equation} \label{dyn2}
\dot Q (t)=\nabla E (\Pi (t)), \qquad \dot\Pi (t)=- \nabla V (Q (t)).
\end {equation}
The main result in \cite{KKS1999} is the following theorem.
Let us denote by $\Vert\cdot\Vert_R$
the norm in $L^2(B_R)$, where $B_R$ is the ball $\{x\in\R^3:|x|<R\}$, and denote $\pi_v(x):=-v\nabla\phi_v(x)$.

\begin {theorem} \label{t7}
Let the condition \eqref{asolV} hold,
and assume that
$
(\psi(0),\,\dot\psi(0),\,q(0),\,p(0))\in\mathcal{S}
$
is a soliton with the total momentum  $ P(0) $. Then the
corresponding solution
$
(\psi(x,t),\,\dot\psi(x,t),\, q (t),\, p (t))
$
of
the system \eqref {wp} with $V=0$ admits the
``adiabatic asymptotics'' 
\[
| q (t) -Q (t) | \le C_0, \quad | P (t) - \Pi (t) | \le C_1 \varepsilon \quad \mbox {\rm for} \quad | t | \le C \varepsilon^{- 1}, 
\]
\[
\sup_{t \in\R} \Big [\Vert \nabla\psi (q (t) + y,t) - \nabla\psi_{v (t)} (y)] \Vert_R
+ \Vert \dot\psi (q (t) + y,t)-\pi_{v (t)} (y) \Vert_R \Big] \le C \varepsilon,
\]
where
$ P (t) $ denotes the total momentum  \eqref {PP}, 
$ v (t)=\mathbf{P}^{- 1} (\Pi (t)) $, and
$(Q (t), \Pi (t)) $ is the solution of the effective Hamiltonian equations  \eqref {dyn2} with initial conditions
$$
Q (0)=q (0), \qquad \Pi (0)=P (0).
$$
\end {theorem}
We note that such relevance of the effective dynamics (\ref{dyn2}) is due to the consistency of Hamiltonian structures:
\medskip \\
I. The effective Hamiltonian  (\ref{Heff}) is the restriction of the Hamiltonian functional (\ref{Ham}) with $ V=0 $
to the solitary manifold $ \mathcal{S} $.
\smallskip \\
II. As shown in \cite{KKS1999}, the
canonical differential form  of the Hamiltonian system\index{Hamiltonian! equation}
(\ref{dyn2}) is also the restriction to $ \mathcal{S} $ of the
canonical differential form
 of the  system \eqref {wp}: formally,
$$
P\cdot dQ=\Bigl [p\cdot dq + \int 
 \pi(x)\,d \psi(x)\,dx \Bigr] \Big |_{\mathcal{S}}.
$$
Therefore, the total momentum $ \mathbf{P} $ is canonically conjugate to the variable $ Q $ on the solitary manifold $ \mathcal{S} $.
This fact justifies the definition (\ref{Heff}) of the effective Hamiltonian as a function of the total
momentum
$ P_v $, and not of the particle momentum $ p_v $.
\medskip

One of the important results in \cite{KKS1999} is the following ``effective dispersion relation'':
\begin {equation} \label{EP}
E (\Pi) \sim \frac {\Pi^2} {2 (1 + m_e)} + {\rm const}, \qquad | \Pi | \ll 1.
\end {equation}
It means that the nonrelativistic mass of a slow soliton increases,
due to the interaction with the field,
by the amount
\begin {equation} \label{me}
m_e=- \frac 13 \langle \rho, \Delta^{- 1} \rho \rangle.
\end {equation}
This increment is proportional to the field energy of the soliton at rest,
\begin{equation}\label{Eown}
E_{\rm own}=H(\Delta^{- 1} \rho,0,0,0)=-\frac 12 \langle \rho, \Delta^{- 1} \rho
\rangle,
\end{equation}
 which
agrees with the Einstein mass-energy equivalence principle (see Section \ref{s4} below).

\begin {remark} \label{rad}
{\rm
The relation (\ref{EP}) suggests that $ m_e $ is an increment of the effective mass.
The true {\it dynamical
justification} for such an interpretation is given by 
Theorem \ref{t7}
which demonstrates the relevance of the effective dynamics (\ref{dyn2}).
}
\end {remark}

\noindent{\bf Generalizations.}
In \cite{KS2000ad}, Kunze and Spohn extended
Theorem \ref{t7} to solitons of the 
Maxwell--Lorentz equations \eqref {mls3} with slowly varying 
potential (\ref{asolV}).
Following the articles
 \cite{KKS1999,KS2000ad}, suitable adiabatic effective dynamics was obtained  in \cite{FTY2002} and  \cite{FGJS2004} for  nonlinear Hartree
and Schr{\"o}dinger equations with slowly varying external potentials, and
in \cite{DS2009, LS2009} and \cite{S2010} for the nonlinear equations of
Einstein--Dirac, Chern--Simon--Schr{\"o}dinger and Klein--Gordon--Maxwell with small external fields. 
In
 \cite{IKS2004b}, the adiabatic  effective dynamics was established for the 
Maxwell--Lorentz equations with a rotating particle.
Similar
adiabatic effective dynamics  was established in \cite{BCFJS} for an electron in
the second-quantized Maxwell field in the presence of a slowly varying external potential.
The results
of
 numerical simulation \cite{KMV2004}
 confirm the adiabatic effective dynamics of solitons for
relativistic 1D 
nonlinear wave equations  (see Section \ref{s9}).

\subsection {Mass-energy equivalence\index{mass-energy equivalence}} \label{s4}
In \cite{KS2000ad}, Kunze and Spohn have established
that in the case
of the Maxwell--Lorentz equations (\ref{mls3})
with slowly varying potential  (\ref{asolV}),
the increment of nonrelativistic mass
also turns out to be proportional to the energy of the static soliton's own field.
Such an equivalence of the self-energy\index{self-energy} of a particle with its mass was first discovered in 1902 by M. Abraham\index{Abraham model}:
he showed by direct calculation that the
energy\index{electrostatic! energy} $ E_{\rm own} $
of electrostatic field\index{electrostatic! field} of an electron at rest adds
\begin{equation}\label{Ame} 
m_e=\displaystyle \frac43 E_{\rm own} / c^2 
\end{equation}
 to its nonrelativistic mass\index{nonrelativistic! mass}
(see \cite{A1902, A1905}, and also \cite[pp.\ 216--217] {K2013}).
By (\ref{Eown}),
this self-energy\index{self-energy} is infinite for a point particle
with the charge density $ \rho(x)=\delta(x) $.
 This means that the field mass\index{field mass} for a point electron\index{point electron} is infinite, which contradicts experiments.
That is why M.~Abraham\index{Abraham model} introduced the model of electrodynamics with
\emph{extended electron}\index{extended electron} \eqref {mls3},
whose self-energy\index{self-energy} is finite.

At the same time, M. Abraham\index{Abraham model} conjectured that the
\emph{entire mass}\index{electron! mass} of an electron is due to its own electromagnetic energy\index{self-energy},
i.e., $ m=m_e $: ``matter disappeared, only energy remains'',
as philosophically-minded contemporaries wrote \cite[pp.~63, 87, 88] {H1908}. (smile :)) 

In 1905, 
the formula (\ref{Ame}) was corrected by A. Einstein, who discovered the famous universal relation
$ E=m_0 c^2 $\index{Einstein!mass-energy equivalence}
which follows from the Special Theory of
Relativity \cite{Einstein-1905}.
Thus, Abraham's discovery of the mass-energy equivalence anticipated  Einstein's 
theory.
The doubtful
 factor $ \frac43 $ in Abraham's\index{Abraham model} formula is due to the nonrelativistic character of the system \eqref {mls3}.
According to the modern view, about 80\% of the electron
mass\index{electron! mass} is of electromagnetic origin \cite{F1966}.

\section{\ Stability of stationary orbits and solitons}

In 1987--1990, Grillakis, Shatah and Strauss developed general theory of orbital stability
of the solitary waves \cite{GSS}.
In \cite{BG}, Bambusi and Galgani established orbital stability 
of solitons for the Maxwell--Lorentz system (\ref{mls3}).

Asymptotic stability of solitary manifold has the meaning
of the local attraction, i.e.,
the convergence to this manifold of states that were sufficiently close to it.
The first results of such type were obtained by  Morawetz,
Segal, and Strauss in the case when the attractor 
consists of the single point zero.
In this case, the attraction is called {\it local energy decay}
\cite{Mor1968}--\cite{St81-2}.

In the case of a continuous attractor,
the key peculiarity of this convergence
 is the instability of the dynamics \emph{along the manifold}.
 The instability  follows directly from the fact that solitons move with different speeds and therefore run away for large times.
Analytically, this instability is caused by the presence of the
Jordan block corresponding to
eigenvalue $ \lambda = 0 $ in the spectrum
of the generator of linearized dynamics. 
The tangent vectors to the soliton manifold are eigenvectors
and generalized eigenvectors corresponding to this Jordan block.
Thus, Lyapunov's theory is not applicable to this case.

\begin{remark}
According to Arnold (see e.g. \cite{AK1998} and references therein),
a similar instability mechanism
may be responsible
for the mixing and ergodicity
of the turbulent flows.
\end{remark}

In a series of articles published during 1985--2003,
Weinstein, Soffer, Buslaev, Perelman, and Sulem
discovered an original  strategy 
for proving asymptotic stability of solitary manifolds.
This strategy relies on i)  special projection of a trajectory onto the solitary manifold,
ii) modulation equations for  parameters of the projection, and iii) time decay of the transversal component.
This approach is a far-reaching development of the Lyapunov stability theory.


\subsection{Spectral stability of solitons}
\label{sect-linear}
Spectral (or linear) stability of solitary waves
can be traced back to 1967,
to the work of Zakharov \cite{zakharov-1967-ak}
on nonlinear Schr\"odinger equation
\begin{eqnarray}\label{nls}
i\dot\psi=-\Delta\psi-|\psi|^2\psi,
\qquad
\psi(x,t)\in\C,
\quad
x\in\R^n,
\quad
n\ge 1,
\end{eqnarray}
which was later developed by Kolokolov \cite{kolokolov-1973-ak}.
This is the ``weakest'' type of stability:
one considers the perturbation of a solitary wave solution
$\phi(x)e^{-i\omega t}$, $\phi(x)\in\R$,
in the form of the Ansatz
\[
\psi(x,t)=(\phi(x)+\chi(x,t))e^{-i\omega t},
\]
with complex-valued
perturbation $\chi(x,t)=u(x,t)+i v(x,t)$
(with $u,\,v$ real-valued),
which is assumed to be small at $t=0$.
Then one
derives the linear approximation to the evolution of $\chi$,
of the form
\begin{eqnarray}\label{l0-l1}
\frac{\partial}{\partial t}
\begin{bmatrix}u\\v\end{bmatrix}
=
A
\begin{bmatrix}u\\v\end{bmatrix},
\end{eqnarray}
and studies the spectrum of the operator
$A=\begin{bmatrix}0&L_0\\-L_1&0\end{bmatrix}$,
where
\[
L_0=-\Delta-\phi^2+\omega,
\qquad
L_1=-\Delta-3\phi^2+\omega
\]
are the Schr\"odinger-type operators.
One can show that
in the case $n\le 2$
one has $\sigma(A)\subset i\R$
(and the corresponding solitary wave is called spectrally stable),
while for $n\ge 3$ the intersection $\sigma(A)\cap(0,+\infty)$
contains a positive eigenvalue
(the corresponding solitary wave is called linearly unstable).
We note that the system \eqref{l0-l1} can not be written as a $\C$-linear
equation on $\chi$:
in the Physics terminology,
the $\mathbf{U}(1)$-symmetry is broken by the choice
$\phi(x)\in\R$.

In the context of the nonlinear Schr\"odinger and Klein--Gordon equations
and similar systems
the spectral stability has been extended to orbital stability;
see
\cite{cazenave1982orbital,GSS,shatah1983stable,W1985}.
At the same time,
the orbital stability for the nonlinear Dirac equation,
and likewise for the Dirac--Maxwell and Dirac--Klein--Gordon systems,
does not seem possible
except via the proof of asymptotic stability,
with the exception
of the completely integrable massive Thirring model
in one spatial dimension \cite{pelinovsky2014orbital-ak}.
Let us mention that the spectral stability for the nonlinear Dirac
equation
in the charge-subcritical cases
has been proved in \cite{BC2019} (see also \cite{BC2019book}).

\subsection {Asymptotic stability of stationary orbits. Orthogonal projection}

This strategy was formed
in 1985--1992 in the pioneering work of Soffer and Weinstein
\cite{SW1990, SW1992,W1985}; see the review \cite{soffer2006}.
The results concern nonlinear $\mathbf{U}(1)$-invariant Schr\"odinger equations with real-valued potential $V(x)$,
\begin {equation} \label {Su}
i \dot \psi(x,t) =-\Delta \psi(x,t) + V(x) \psi(x,t) + \lambda | \psi(x,t) |^p \psi(x,t), \qquad x \in \R^n,
\end {equation}
where $ \lambda \in \R $, $ p = 3 $ or $ 4 $, $ n = 2 $ or $ n = 3 $,
and $ \psi(x,t) \in \C $.
The corresponding Hamiltonian functional is given by
$$
\mathcal{H} = \int \Big[\frac12 | \nabla \psi(x) |^2 + \frac12V(x) | \psi(x) |^2 + \frac \lambda p | \psi(x) |^p\Big]\,dx.
$$
For $ \lambda = 0 $, equation (\ref {Su}) is linear. Let $ \phi_*(x) $ denote its ground state corresponding to the minimal eigenvalue
$ \omega_* <0 $. Then $ C \phi_*(x) e^{- i \omega_* t} $
are periodic solutions for any complex constant $C$.
Corresponding phase curves are circles filling the complex plane.
For nonlinear equations (\ref {Su}) with a small real $ \lambda \ne 0 $, it turns out that a wonderful
{\it bifurcation} occurs:
small neighborhood of zero in the complex plane
turns into an analytic
soliton manifold
$ \mathcal{S} $
invariant with respect to the dynamics group corresponding to \eqref{Su},
which is still filled with invariant
circles
$ \phi_\omega(x) e^{- i \omega t} $
whose frequencies $ \omega $ are close to $ \omega_* $.

The  main result of \cite{SW1990, SW1992} (see also \cite{PW1997}) is  long-time attraction to one of these circles
 for any solution with sufficiently small initial data:
\begin {equation} \label {soli}
\psi(x,t) = \phi_{\pm}(x) e^{- i \omega_{\pm} t} + r_\pm(x,t),
\end {equation}
where the remainder decay in weighted norms: for any $ \sigma> 2 $,
$$
\Vert \langle x \rangle^{- \sigma} r_\pm (\cdot, t) \Vert_{L^2 (\R^n)} \to 0, \qquad t \to \pm \infty,
$$
where $ \langle x \rangle: = (1+ | x |)^{1/2} $.
The proof relies on  linearization of the dynamics and decomposition of solutions into two components,
$$
\psi (t) = e^{- i \Theta (t)} (\phi_{\omega (t)} + \chi (t)),
$$
satisfying the following orthogonality condition \cite[(3.2) and (3.4)] {SW1990}:
\begin {equation} \label {or}
\langle \phi_{\omega (t)}, \chi (t) \rangle = 0.
\end {equation}
This orthogonality and dynamics (\ref {Su}) imply the {\it modulation equations} for $ \omega (t) $ and $ \gamma (t) $, where
$ \gamma (t): = \Theta (t)-\displaystyle \int_0^t \omega (s) ds $
(see 
\cite[(3.9a)--(3.9b)]{SW1990}).
The orthogonality (\ref {or}) implies that the component $ \phi (t) $ lies in the continuous spectral space of the Schr\"odinger operator
\[
H (\omega_0): =-\Delta + V + \lambda |\phi_{\omega_0}|^p,
\]
which leads to the time decay of  $\chi (t)$ (see \cite[(4.2a) and (4.2b)] {SW1990}).
Finally, this decay implies the convergence $ \omega (t) \to \omega_\pm $ and the asymptotics (\ref {soli}).
\smallskip

These results and methods were subsequently widely developed
for the nonlinear  Schr\"odinger, wave and Klein--Gordon equations with potentials under various spectral assumptions on linearized dynamics;
see, for instance, \cite{SW1999,SW2004}.

\medskip

\subsection {Asymptotic stability of solitons. Symplectic projection}
\begin {figure} [htbp]
\begin {center}
\includegraphics [width = 0.8 \columnwidth] {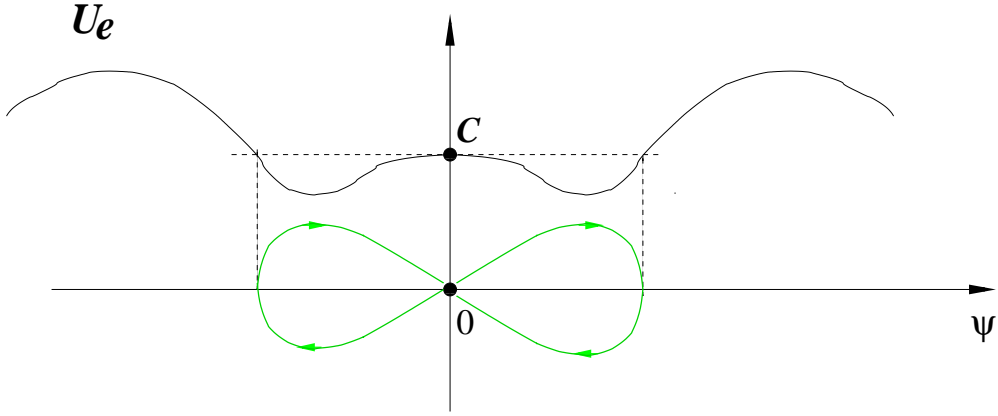}
\caption {Reduced potential and soliton.}
\label {Ue}
\end {center}
\end {figure}

Next breakthrough  in the theory of asymptotic stability was achieved
in 1990--2003 
by Buslaev, Perelman and Sulem \cite{BP1993}--\cite{BS2003},
who first proved  asymptotics of type (\ref {soli}) for translation-invariant 1D Schr\"odinger equations 
\begin {equation} \label {BPS}
i \dot \psi(x,t) =-\psi ''(x,t) -F (\psi(x,t)), \qquad x \in \R,
\end {equation}
which are also assumed to be $\mathbf{U}(1) $-invariant, i.e.,
$F(\psi)=-\nabla U(\psi)$, and
\begin{equation}\label{C3}
F(e^{i\theta}\psi)=e^{i\theta}F(\psi), \qquad \psi\in\C,\,\,\,\theta\in\R.
\end{equation}
Moreover, they assume that
\begin{equation}\label{F10}
U(\psi)=\mathcal{O}(|\psi|^{10}),\qquad \psi\to 0.
\end{equation}
 Under some simple conditions on the potential $U$, 
there exist finite energy solutions of  the form
\begin{equation}\label{sol0}
 \psi(x,t) = \phi_0(x) e^{ i \omega_0 t},
 \end{equation}
 with $ \omega_0>0 $. The amplitude $ \phi_0(x) $ satisfies the corresponding stationary equation
 \begin{equation}\label{cse}
 -\omega_0 \phi_0(x) =-\phi_0 ''(x) -F (\phi_0(x)), \qquad x \in \R,
 \end{equation}
which implies the ``conservation law''
 \begin{equation}\label{solC}
 \frac {| \phi_0 '(x) |^2} 2 + U_e (\phi_0(x)) = E,
 \end{equation}
where the effective potential is given by
$U_e (\psi) = U(\psi) + \omega_0 \frac {| \psi |^2}2 \sim \omega_0 \frac {| \psi |^2}2$ as  $\psi\to 0$
by (\ref{F10}).
 For the existence of a finite energy solution,  the graph of ``effective potential'' $ U_e (\psi) $ 
should be similar to Fig. \ref {Ue}. 
The finite energy solution
$\phi_0(x)$
 is defined by (\ref {solC})  with the constant $ E = U_e (0) $
 since for other $E$ the solutions to
 (\ref{solC})
 do not converge to zero as $|x|\to\infty$. 
 Equation (\ref{solC}) with  $ E = U_e (0) $ implies that
 \begin{equation}\label{solC2}
 \frac {| \phi_0 '(x) |^2} 2 = U_e (0) -U_e (\phi_0(x))
\sim \frac{\omega_0}  2\phi_0^2(x).
 \end{equation}
 Hence, for finite energy solutions, one has
 \begin{equation}\label{rU}
 \phi_0(x) \sim e^{- \sqrt{\omega_0} | x |},\qquad |x|\to\infty.
 \end{equation}

It is easy to verify that  the  following functions are also solutions  ({\it moving solitons}):
\begin{equation}\label{solv}
 \psi_{\omega, v, a, \theta}(x,t) = \phi_\omega (x-vt-a) e^{ i (\omega t + kx + \theta)}, \qquad
 \omega = \omega_0-v^2/4, \qquad k = v / 2.
 \end{equation}
The set of all such solitons with parameters $ \omega, v, a, \theta $ forms a 4-dimensional smooth submanifold $ \mathcal{S} $ in the
Hilbert phase space $ \mathcal{X}: = L^2 (\R) $. Solitons (\ref {solv}) are obtained from  (\ref {sol0}) by the Galilean transformation
\begin{equation}\label{Gavt}
G (a, v, \theta): \psi(x,t) \mapsto \varphi(x,t) =\psi (x-vt-a, t) e^{ i (-\frac {v^2} 4t +\frac v2 x + \theta)};
\end{equation}
the Schr\"odinger equation (\ref {BPS}) is invariant with respect to this group of transformations.

Linearization of the Schr\"odinger equation (\ref {BPS}) at the soliton (\ref {sol0}) is obtained as in Section~\ref{sect-linear},
by substitution
$ \psi(x,t) = (\phi_0(x) + \chi(x,t)) e^{i \omega_0t} $ and retaining terms of the first order in $ \chi $.
The resulting equation
is not linear over the field of complex numbers, since it contains
both $ \chi $ and $ \overline \chi $. 
This follows from the fact that the nonlinearity
$ F (\psi)=a(|\psi|^2)\psi $ is not complex-analytic due to
 the $\mathbf{U}(1)$-invariance 
(\ref {C3}).  Complexification of this linearized equation reads
as follows:
\begin{equation}\label{cLS}
 \dot \Psi(x,t) = A_0\Psi(x,t),\qquad A_0=-jH_0,
\end{equation}
where $ j $ is a real $ 2 \times 2 $ matrix, representing the multiplier $ i $, $ \Psi(x,t) \in \C^2 $,
and  $ H_0 =-{d^2} / {dx^2} +\omega_0 + V(x) $, where $ V(x) $
is a real, matrix-valued potential
which decreases exponentially as $ | x | \to \infty $ due to (\ref {rU}). 
Note that the operator $A_0 $
in \eqref{cLS}
corresponds to the linearization   at the soliton
(\ref{solv}) with parameters $\omega=\omega_0$ and $a=v=\theta=0$.
Similar operators $A(\omega, a, v, \theta) $  corresponding to linearization at solitons (\ref{solv}) with various parameters $ \omega, a, v, \theta $ 
are also connected by the linear Galilean transformation (\ref {Gavt}). 
Therefore, their spectral properties completely coincide. In particular, their continuous spectra coincide with $(-i\infty,-i \omega_0]\cup [i\omega_0, i\infty) $.

Main results of \cite{BP1993}--\cite{BS2003} are asymptotics of type (\ref {soli})
for solutions with initial data close to the solitary manifold $ \mathcal{S} $:
\begin {equation} \label {sdw}
\psi(x,t) = \phi_\pm (x-v_\pm t) e^{- i (\omega_\pm t + k_\pm x)} + W (t) \Phi_\pm + r_\pm(x,t),\qquad \pm t>0,
\end {equation}
where $ W (t) $ is the dynamical group of the free Schr\"odinger equation, $ \Phi_{\pm} $ are some scattering states of finite energy, and
$ r_{\pm} $ are  remainder  terms which decay to zero in the global norm:
\begin {equation} \label {glob2u}
\Vert r_{\pm} (\cdot, t) \Vert_{L^2 (\R)} \to 0, \qquad t \to \pm \infty.
\end {equation}
These asymptotics were obtained under the following assumptions on the spectrum of the generator $B_0$:
\smallskip

U1. The discrete spectrum of the operator $ C_0 $ consists of exactly three eigenvalues $ 0 $ and $ \pm i\lambda 0 $, and
\begin{equation}\label{fot}
 \lambda< \omega_0<2 \lambda.
\end{equation}
This condition means that the discrete mode interacts with 
the continuous spectrum already  in the first order of approximation.

U2. The edge points $\pm i\omega_0$ of the continuous spectrum  are neither  eigenvalues, nor virtual levels (threshold resonances) of $ C_0 $.

U3. Furthermore, it is assumed the condition \cite[(1.0.12)] {BS2003}, which means a strong coupling of discrete and continuous spectral components, 
providing energy radiation, similarly to the Wiener condition  \cite[(6.24)]{K2021}, \cite[(1.5.13)]{KK2021}.
It ensures that the interaction of  discrete component with  continuous spectrum does not vanish already in the first order of perturbation theory.
 
This condition is a nonlinear version of the Fermi Golden Rule \cite{RS4}, which was introduced by Sigal 
in the context of nonlinear PDEs \cite{Sig1993}.
\smallskip

Examples of potentials satisfying all these conditions are constructed in \cite{KKK2013}.  
\smallskip

In 2001, Cuccagna extended results of  \cite{BP1993}--\cite{BS2003} to translation-invariant Schr\"odinger equations 
in $\R^n$, $ n \ge 2 $ \cite{C2001}.
Similar approach 
to the asymptotic stability of solitons
was developed
for the Korteweg--de Vries equation in \cite{PW1994},
for the regularized long-wave equation
in \cite{MW1996},
and
for the nonlinear Schr\"odinger equation in $\R^3$ with concentrated nonlinearity
in \cite{ANO}.
In the context of the nonlinear Dirac equation with a potential,
stability with respect to perturbations
of solitary waves in particular directions
was considered in \cite{B2006}.


\medskip

\noindent
{\bf Method of symplectic projection in the Hilbert  phase space.}
Novel approach \cite{BP1993}--\cite{BS2003} relies on {\it symplectic projection} $ P $ of solutions onto  the solitary manifold. This means that 
$$
Z: = \psi-S \quad \mbox {symplectic-orthogonal to the tangent space} \quad \mathcal{T}: = T_{S} \mathcal{S}
$$
for the projection $ S: = P \psi $. This projection is correctly defined in a small neighborhood of $ \mathcal{S} $. 
For this it is important that $ \mathcal{S} $ is a {\it symplectic manifold}, 
i.e. the corresponding symplectic form is non-degenerate on the tangent spaces $T_{S}\mathcal{S}$.

In particular, 
the  approach \cite{BP1993}--\cite{BS2003} allowed to get rid of  the assumption that the initial data are small.

Now a solution $ \psi(t)$ for each $ t> 0 $ decomposes onto {\it symplectic-orthogonal} components  $\psi(t)= S(t)+ Z(t)$, 
where $ S (t): = P \psi (t) $, and the dynamics is linearized on the soliton $ S (t)$.  

The Hilbert phase space $ \mathcal{X}: = L^2 (\R) $ is split as $ \mathcal{X} = \mathcal{T} (t) \oplus \mathcal{Z} (t) $, where
$ {\mathcal Z} (t) $ is {\bf symplectic-orthogonal} complement to the tangent space $ \mathcal{T} (t): = T_{S (t)} \mathcal{S} $.
The corresponding equation for the {\it transversal component} $ Z (t) $ reads
$$
\dot Z (t) = A (t) Z (t) + N (t),
$$
where $ A (t) Z (t) $ is the linear part, and $ N (t) = \mathcal{O} (\Vert Z (t) \Vert^2) $ is the corresponding nonlinear part.
The main difficulties in studying this equation are as follows i) it is {\it non-autonomous},
and ii) the generators $ A (t) $ {\it  are not self-adjoint} (see Appendix in \cite{KKsp2014}).
It is important that $ A (t) $ are {\it Hamiltonian operators}, for which the existence of spectral  decomposition 
is provided by the Krein--Langer theory of $ J $ -selfadjoint operators \cite{KL1963, L1981}.
In \cite{KKsp2014, KKsp2015} we have developed a special version of this theory  providing
the corresponding eigenfunction expansion which is necessary for the justification 
 of the  approach \cite{BP1993}--\cite{BS2003}.
\smallskip \\
The main steps of the  strategy \cite{BP1993}--\cite{BS2003} are as follows:

\smallskip

\noindent
$\bullet$ {\bf Modulation equations.}
The parameters of the soliton $ S (t) $ satisfy  {\bf modulation equations}: 
for example, for the  speed  $v (t)$ we have
$$ 
\dot v (t) = M (\psi (t)), 
$$ 
where $ M (\psi) = \mathcal{O} (\Vert Z \Vert^2) $ for small norms $ \Vert Z \Vert $. 
Therefore, the parameters change ``super-slow'' near the soliton manifold, like adiabatic invariants.

\smallskip

\noindent
$\bullet$ {\bf Tangent and transversal components.}
The transversal component $ Z (t) $ in the splitting $ \psi (t) = S (t) + Z (t) $ belongs to the transversal space $ \mathcal{Z} (t) $.
The tangent space $ \mathcal{T} (t) $ is the root space of the generator $ A (t) $ 
and corresponds to the ``unstable'' spectral point $ \lambda = 0 $. The key observation is that

i) symplectic-orthogonal subspace
$ {\mathcal Z} (t) $
does not contain ``unstable'' tangent vectors, and, moreover,

ii) the subspace $ {\mathcal Z} (t) $ is {\bf invariant} with respect to the generator $ A (t) $, since
the subspace $ \mathcal{T} (t) $ is invariant, and $ A (t) $  is the  Hamiltonian operator.

\smallskip

\noindent
$\bullet$ {\bf Continuous and discrete components.}
The transversal component allows further splitting $ Z (t) = z (t) + f (t) $, where $ z (t) $ and $ f (t) $ 
belong, respectively, to discrete and continuous spectral subspaces $ \mathcal{Z}_d (t) $ and $ \mathcal{Z}_c (t) $ of $A(t)$
in the space $ \mathcal{Z} (t) = \mathcal{Z}_d (t) + \mathcal{Z}_c (t) $.

\smallskip

\noindent
$\bullet$ {\bf  Poincare normal forms and  Fermi Golden Rule.}
The component $z (t) $ satisfies a nonlinear equation, which is reduced to Poincare normal form 
up to higher order terms \cite[Equations (4.3.20)]{BS2003}.
(For the relativistic-invariant Ginzburg--Landau equation, a similar reduction done in \cite[Equations (5.18)]{KopK2011b}.)  
The normal form allowed  to derive some ``conditional'' decay for $z(t) $ using the Fermi Golden Rule
\cite[(1.0.12)] {BS2003}.

\smallskip

\noindent
$\bullet$ {\bf  Method of majorants.}
A skillful combination of the conditional decay for $ z (t) $ with super-slow evolution of the soliton parameters
allows us to prove the decay for $ f (t) $ and $ z (t) $ by the method of  majorants.
Finally, this decay implies the asymptotics (\ref {sdw})--(\ref {glob2u}).

\subsection {Generalizations and Applications} \label{s8}

\noindent
{\bf $N$-soliton solutions.} The methods and results of \cite{BS2003} were developed in \cite{M2005}--\cite{MMT2002} and  \cite{P2004}--\cite{RSS2005}
for $ N $-soliton solutions for translation-invariant nonlinear Schr\"odinger equations.

\smallskip

\noindent
{\bf  Multiphoton  radiation.}
In  \cite{CM2008}, Cuccagna and Mizumachi  extended methods  and results  of \cite{BS2003}  to the case 
when the inequality (\ref {fot}) is changed to
$$
N  \lambda< \omega_0< (N+1) \lambda,
$$
with some natural $N>1$,  and the corresponding analogue of condition U3 holds.
It means that the interaction of discrete modes with a continuous spectrum occurs only in the $ N $-th order of 
perturbation theory. The decay rate of the remainder term (\ref {glob2u}) worsens
 with growing  $ N $.

\smallskip

\noindent
{\bf Linear equations coupled to nonlinear oscillators and particles.}
In \cite{BKKS2008, KKopSt2011} methods and results of \cite{BS2003} were extended
to i) the Schr\"odinger equation interacting with a nonlinear
$\mathbf{U}(1) $-- invariant oscillator, 
ii) in \cite{IKS2011, IKV2011}-- to 
translation-invariant systems of  relativistic particle coupled to
the wave and Maxwell equations, 
and iii) in \cite{IKV2006, KKop2006, KKopS2011}-- to  similar translation-invariant systems of  relativistic particle coupled to
the Klein--Gordon, Schr\"odinger and Dirac equations.
The review of the results \cite{IKV2006}--\cite{IKV2011} can be found in \cite{Im2013}.

\smallskip

\noindent
{\bf Relativistic equations.}
In \cite{K2002, BC2012, Kumn2013, KopK2011a, KopK2011b} methods and results \cite{BS2003} were extended for the first time 
to {\it relativistic-invariant} nonlinear equations.  Namely, i) in \cite{K2002} and \cite{Kumn2013}--\cite{KopK2011b} asymptotics of type (\ref {sdw}) 
were obtained for 1D relativistic-invariant nonlinear wave equations with potentials of the Ginzburg--Landau type,
\begin{equation}\label{GL1}
\ddot\psi(x,t)=\Delta \psi(x,t)-m^2\psi(x,t)+f(\psi(x,t)),
\end{equation}
and ii)  in \cite{BC2012} and \cite{CPS2017}
for relativistic-invariant nonlinear Dirac equations.
In  \cite{KKK2013} we have constructed examples of potentials providing all spectral properties of the linearized dynamics
imposed in \cite{Kumn2013}--\cite{KopK2011b}.

\smallskip

\noindent
{\bf The justification of the eigenfunction expansions.}
In \cite{KKsp2014,KKsp2015} we have justified the eigenfunction expansions  for nonselfadjoint Hamiltonian operators 
which were used in \cite{Kumn2013}--\cite{KopK2011b}. For the justification we have developed a
special version of the Krein--Langer  theory of $ J $ -selfadjoint operators \cite{KL1963, L1981}

\medskip

\noindent
{\bf Vavilov--Cherenkov radiation.}
The article \cite{FG2014} concerns the nonrelativistic particle  
coupled to
the Schr\"odinger equation  (system (1.9)--(1.10) in \cite{FG2014}). 
This system  is considered as a model of the Cherenkov  radiation. 
The main result of \cite{FG2014} is long-time convergence to a soliton with the sonic speed for initial solitons with a supersonic speed 
in the case of a weak interaction (``Bogoliubov limit'') and small initial field.

Asymptotic stability of solitons for  similar system was established in \cite{KKop2006}.

\subsection {Further generalizations}
The results on asymptotic stability of solitons were developed in different directions.

\smallskip

\noindent
{\bf Systems with several bound states.}
The case of many simple eigenvalues of linearization at a soliton (\ref {cLS}) was first investigated in 
\cite{BC2011}--\cite{TY2002} for  the nonlinear Schr\"odinger equation
$$
i \dot \psi(x,t) = (- \Delta + V(x)) \psi(x,t) \pm | \psi(x,t) |^2 \psi(x,t), \qquad x \in \R^3.
$$
Asymptotic stability and long-time  asymptotics of solutions were established under the following assumptions:
\smallskip

i) the endpoint of the essential spectrum
is neither an eigenvalue nor a
virtual level
(threshold resonance) for linearized equation;

ii) the eigenvalues of the linearized equation satisfy the corresponding non-resonance condition;

iii) a new version of the Fermi Golden Rule.

\smallskip

The main  obtained result: any solution with small initial data converges to some ground state at $ t \to \infty $
with speed $ t^{- 1/2} $ in $ L^2_{\rm loc} (\R^3) $.
There are different long-time regimes depending on relative contributions of eigenfunctions into initial data.

One of the difficulties is the possible existence of invariant tori corresponding to the resonances. 
Great efforts have been  applied to show that the role of  metastable tori decays like
$t^{- 1/2}$ as $ t \to \infty $.

This result was extended in \cite{BC2011} to nonlinear Klein--Gordon  equations
$$
\ddot \psi(x,t) = (\Delta-V(x) -m^2) \psi(x,t) + f(\psi(x,t)), \qquad x \in \R^3.
$$
Namely, under conditions i) -- iii) above, any sufficiently small solution for large times asymptotically is a free wave 
in the norm $H^1(\R^3)$. The proofs  largely rely on the theory of the Birkhoff  normal forms. The main novelty is 
the transition to normal forms without  loss of the  Hamiltonian structure.

In \cite{C2011} the results of \cite{CM2008, BC2011} are extended to nonlinear Schr\"odinger equation
$$
i \dot \psi(x,t) = (- \Delta + V(x)) \psi(x,t) + \beta (| \psi(x,t) |) \psi(x,t), \qquad x \in \R^3.
$$
Main results are  long-time asymptotics ``ground state + dispersive wave''
in the norm $ H^1 (\R^3) $ for solutions  close to the ground state.

The corresponding linearized equation can have many eigenvalues, which must satisfy 
the non-resonant conditions \cite{BC2011}, and  suitable modification of the Fermi Golden Rule holds.
However, for nonlinear Schr\"odin\-ger equations  methods  of \cite{BC2011} required a significant improvement: now
the canonical coordinates are constructed using the Darboux theorem.

\smallskip

\noindent
{\bf General theory of relativity.}
The article \cite{HM2004} concerns so-called ``kink instability'' of  self-similar and spherically symmetric
solutions of the equations of the general theory of relativity with a scalar field, as well as with a ``hard fluid'' as sources. 
The authors constructed examples of self-similar solutions that are unstable to the kink perturbations.

The article \cite{DR2011} examines linear stability of slowly rotating Kerr solutions for the Einstein equations in vacuum.
In \cite{T2013} a pointwise damping of solutions to the wave equation is investigated for the case 
of stationary asymptotically flat space-time in the three-dimensional case.

In \cite{AB2015} the Maxwell equations are considered outside  slowly rotating Kerr black hole.
The main results are: i) ~ boundedness of a positive definite energy on each hypersurface $ t = \const $ 
and ii) convergence of each solution to a stationary Coulomb field.

In \cite{DSS2012} the pointwise decay was proved for linear waves against the Schwarzschild black hole.

\smallskip

\noindent
{\bf Method of concentration compactness.}
In \cite{KM2006} the concentration compactness method was used for the first time to prove global well-posedness, 
scattering and blow-up of solutions to critical focusing nonlinear Schr\"odinger equation
$$
i \dot \psi(x,t) =-\Delta \psi(x,t)-| \psi(x,t) |^{\frac4 {n-2}} \psi(x,t), \qquad x \in \R^n
$$
in the radial case. Later on, these methods were extended in \cite{KM2012, DKM2016, KM2008, KNS2015}
to general non-radial solutions and to nonlinear wave equations of the form
$$
\ddot \psi(x,t) = \Delta \psi(x,t) + | \psi(x,t) |^{\frac4 {n-2}} \psi(x,t), \qquad x \in \R^n.
$$
One of the main results is splitting of the set of initial states, close to the critical energy level,
into three subsets with certain long-term asymptotics:
either a blow-up in a finite time, or an asymptotically free wave,
or the sum of the ground state and an asymptotically free wave.
All three alternatives are possible; all nine combinations with $ t \to \pm \infty $ are also possible. 
Lectures \cite{NS2011} give excellent introduction to this area.
The articles \cite{DKM2014, KM2011} concern super-critical nonlinear wave equations.

Recently, these methods and results were extended to critical wave mappings  \cite{KLS2014}--\cite{KS2012}. 
The ``decay onto solitons'' is proved: every $1$-equivariant finite-energy wave mapping
of exterior of a ball with Dirichlet boundary conditions into  three-dimensional sphere
exists globally in time and dissipates into a single stationary solution of its own topological class.

\section{\ Linear dispersion}

In all results on long-time asymptotic  for nonlinear Hamiltonian PDEs,
the key role 
 is played by
dispersive decay of solutions of the corresponding linearized equations.
In this section,
we choose most important or recent articles
out of a huge number of publications.

\medskip

\noindent
{\bf Dispersion decay in weighted Sobolev norms.} 
Dispersion decay for wave equations was first proved in  linear scattering theory \cite{LMP1963}.
For the Schr\"odinger equation with a potential,
a systematic approach to dispersive decay was proposed by 
Agmon \cite{Agmon} and Jensen and Kato \cite{JK}. This theory was extended by many authors
to  wave,  Klein--Gordon and Dirac equations and to the corresponding discrete equations, see
\cite{BG2012}--\cite{EKT}, \cite{JSS1991}--\cite{KT2016}  and references therein. 
\medskip

\noindent
{\it $L^1$-$L^\infty $ decay}
for solutions of the linear Schr\"odinger equation
was proved for the first time by Journ\'e, Soffer, and Sogge \cite{JSS1991}.
Namely, it was shown that the solutions to
\begin {equation} \label {LSE}
i \dot \psi(x,t) = H \psi(x,t): = (- \Delta + V(x)) \psi(x,t), \qquad x \in \R^n,
\qquad
n\ge 3
\end {equation}
satisfy
\begin {equation} \label {l1i}
\Vert P_c \psi (t) \Vert_{L^\infty (\R^n)} \le Ct^{- n / 2} \Vert \psi (0) \Vert_{L^1 (\R^n)}, \qquad t>0,
\end {equation}
provided that $ \lambda = 0 $ is neither an eigenvalue
nor a virtual level of $ H $.
Here $ P_c $ is an orthogonal projection onto continuous  spectral space of the operator $ H $.
It is assumed that the potential $V(x) $
is real-valued,
sufficiently smooth, and rapidly decays as $ | x | \to \infty $.
This result was  generalized later by many authors;
see below.
\smallskip

The decay of type  (\ref {l1i})  and Strichartz estimates were established
in \cite{RS2004}
for 3D Schr\"odinger equations  (\ref {LSE}) with ``rough'' and time-dependent potentials $ V = V(x,t) $.
Similar estimates were obtained in \cite{BG2012} for 3D Schr\"odinger and  wave equations with (stationary) Kato class potentials.
\smallskip

In \cite{EGG2014}, the 4D Schr\"odinger equations (\ref {LSE}) are considered
in the case when there is a
virtual level or an eigenvalue at  zero energy.
In particular, in the case of  an eigenvalue at  zero energy,  there is a time-dependent  operator $ F_t $ of rank $ 1 $
such that $ \Vert F_t \Vert_{L^1 \to L^\infty} \le 1 / \log t $ for $ t> 2 $, and
$$
\Vert e^{itH} P_c-F_t \Vert_{L^1 \to L^\infty} \le Ct^{- 1}, \qquad t> 2.
$$
Similar dispersive estimates were also proved for solutions to 4D wave equation with a potential.
\smallskip

In \cite{GG2015,GG2017}, the Schr\"odinger equation  (\ref {LSE}) is considered in $ \R^n $ with  $ n \ge 5 $
with the assumption that zero is an eigenvalue
(the Schr\"odinger operator in dimension $n\ge 5$
could only have zero eigenvalues
but no genuine virtual level at zero).
It is shown, in particular, that
there is a time-dependent rank one operator $ F_t $  such that
$ \Vert F_t \Vert_{L^1 \to L^\infty} \le C | t |^{2-n / 2} $ for $ | t |> 1 $, and
$$
\Vert e^{itH} P_c-F_t \Vert_{L^1 \to L^\infty} \le C | t |^{1-n / 2}, \qquad | t |> 1.
$$
With a stronger decay of the potential, the evolution admits an operator-valued expansion
$$
e^{itH} P_c (H) = | t |^{2-n / 2} A_{- 2} + | t |^{1-n / 2} A_{- 1} + | t |^{- n / 2} A_0,
$$
where $ A_{- 2} $ and $ A_{- 1} $ are  finite rank operators $ L^1 (\R^n) \to L^\infty (\R^n) $, while $ A_0 $
maps weighted  $ L^1 $ spaces to  weighted  $ L^\infty $ spaces.
The operators $ A_{- 2} $ and $ A_{- 1} $
are equal to zero under certain conditions of the orthogonality of the potential $ V $ to the eigenfunction which corresponds to zero eigenvalue.
Under the same orthogonality conditions, the remainder term $ | t |^{- n / 2} A_0 $ also maps $ L^1 (\R^n) $ to $ L^\infty (\R^n) $, 
and therefore the group $ e^{itH} P_c (H) $ has the same dispersion decay as
the free evolution, despite its eigenvalue at zero.

\medskip

\noindent
{\it $ L^p$-$L^q $ decay}  was first established in \cite{MWS1980} for solutions
of the free Klein--Gordon equation $ \ddot \psi = \Delta \psi- \psi $ with initial state $ \psi (0) = 0 $:
\begin {equation} \label {LpLqKG}
\Vert \psi (t) \Vert_{L^q} \le Ct^{-d} \Vert \dot \psi (0) \Vert_{L^p}, \qquad t> 1,
\end {equation}
where $ 1 \le p \le 2 $, $ 1 / p + 1 / q = 1 $, and $ d \ge 0 $ is a piecewise-linear function of $ (1 / p, 1 / q) $.
The proofs use the Riesz interpolation theorem.
\smallskip

In \cite{BS1993}, the estimates (\ref {LpLqKG}) were extended to solutions of the perturbed Klein--Gordon equation
$$
\ddot \psi = \Delta \psi- \psi + V(x) \psi
$$
with $ \psi (0) = 0 $. It is shown that (\ref {LpLqKG}) holds for $ 0 \le 1 / p-1/2 \le 1 / (n + 1) $. 
The smallest value of $ p $ and the fastest decay rate $ d $ occurs when $ 1 / p = 1/2 + 1 / (n + 1) $, $ d = (n-1) / (n + 1) $. 
The result is proved under the assumption that the potential $ V $ is  smooth and small in a suitable sense. 
For example, the result holds true when $ | V(x) | \le c (1+ | x |^2)^{- \sigma} $, where $ c>0 $ is sufficiently small.
Here  $ \sigma> 2 $ for $ n = 3 $,  $ \sigma> n / 2 $ for  odd $ n\ge 5 $, and $ \sigma> (2n^2 + 3n + 3) / 4 (n + 1) $ for even $ n \ge 4 $. 
The results also apply to the case when $ \psi (0) \ne 0 $.
\smallskip

The seminal article \cite{JSS1991} concerns $ L^p$-$L^q $ -decay of solutions to the Schr\"odinger equation (\ref {LSE}).
It is assumed that $ (1+ | x |^2)^\alpha V(x) $ is a multiplier in the Sobolev spaces $ H^\eta $ for some $ \eta> 0 $ and $ \alpha> n + 4 $,
and the Fourier transform of $ V $ belongs to $ L^1 (\R^n)$. 
Under this conditions,  the main result of \cite{JSS1991} is the following theorem: if $\lambda = 0$ is neither an eigenvalue nor a virtual level of $ H$, 
then
\begin {equation} \label {LpLq}
\Vert P_c \psi (t) \Vert_{L^q} \le Ct^{- n (1 / p-1/2)} \Vert \psi (0) \Vert_{L^p}, \qquad t> 1,
\end {equation}
where $ 1 \le p \le 2 $ and $ 1 / p + 1 / q = 1 $.
Proofs are based on $ L^1$-$L^\infty $-decay (\ref {l1i}) and the Riesz interpolation theorem.

In \cite{Y2005}, under suitable conditions on the spatial decay of $ V(x) $,
the estimates (\ref {LpLq}) were proved for  all $ 1 \le p \le 2 $
if $ \lambda=0 $
is  neither an eigenvalue nor a virtual level of $ H $,
and for all $ 3/2 <p \le 2 $ otherwise.

\medskip

\noindent
{\it The Strichartz estimates} were extended
in \cite{AFVV2010} to the Schr\"odinger equations with the magnetic field in $\R^n$ with $ n \ge 3 $;
in \cite{A2015} -- to the wave equations with a magnetic potential in $\R^n $ for $ n \ge 3 $;
in \cite{BG2014} -- to the wave equation in $\R^3$
with potentials of the  Kato class.


\section{\ Numerical simulation of soliton asymptotics}
\label{s9}
Here we describe the results of Arkady Vinnichenko (1945-2009) on numerical simulation of    i) global 
attraction to solitons \eqref {att} and \eqref {attN},
and ii) adiabatic effective dynamics of solitons  for relativistic-invariant one-dimensional
nonlinear wave equations \cite{KMV2004}.

\subsection {Kinks of relativistic-invariant Ginzburg--Landau equations}
First, we simulated numerically solutions to relativistic-invariant 1D nonlinear wave equations with polynomial nonlinearity
\begin {equation} \label {GL}
\ddot \psi(x,t) = \psi ''(x,t) + F (\psi(x,t)), 
\end {equation}
where $F (\psi) =-\psi^3 + \psi$.
Since $F (\psi) = 0$ for $\psi = 0,\,\pm 1$, there are three  equilibrium states: $S(x)\equiv  0,\, + 1,\, -1$.

The corresponding  potential $ U(\psi) = \frac {\psi^4} {4}-\frac {\psi^2} {2} +\frac 14$ has minima at $ \psi = \pm 1 $ 
and a local maximum at $ \psi = 0 $, therefore two finite energy solutions  $\psi={\pm 1}$ are stable, and the solution $\psi=0$ with infinite energy  is unstable.
Such potentials with two wells are called Ginzburg--Landau potentials.

Besides the constant stationary solutions $ S(x) \equiv 0, + 1, -1 $, there is also a non-constant one,
$ S(x) = \tanh \frac x{\sqrt2} $,
called a ``kink''.
Its shifts and reflections $ \pm S (\pm x-a) $ are also solutions,
as well as
 their Lorentz transforms
$$ 
\pm S (\gamma (\pm x-a-vt)),\qquad  \gamma = 1 / {\sqrt {1-v^2}}, \quad | v | <1. 
$$
These are uniformly moving ``traveling waves'' (i.e. solitons). The kink is strongly compressed when the velocity $ v $ is close to $ \pm 1 $.
Equation (\ref {GL}) is formally equivalent to the Hamiltonian system with the Hamiltonian
\begin {equation} \label {hamGL}
\mathcal{H} (\psi, \pi) = \int 
\Big[\frac12 | \pi(x) |^2 + \frac12 | \psi '(x) |^2 + U(\psi(x))\Big]\,dx.
\end {equation}
This Hamiltonian is finite  for functions $ (\psi, \pi)$ from the Hilbert
space  $\mathcal{E}=H^1(\R^3)\oplus L^2(\R^3) $, 
for which the convergence
$$
\psi(x) \to \pm 1 \quad\mbox{as}\quad | x | \to \pm \infty
$$
is sufficiently fast.

\medskip

\begin{figure}[htbp]
\begin{center}
\includegraphics[width=1.00\columnwidth]{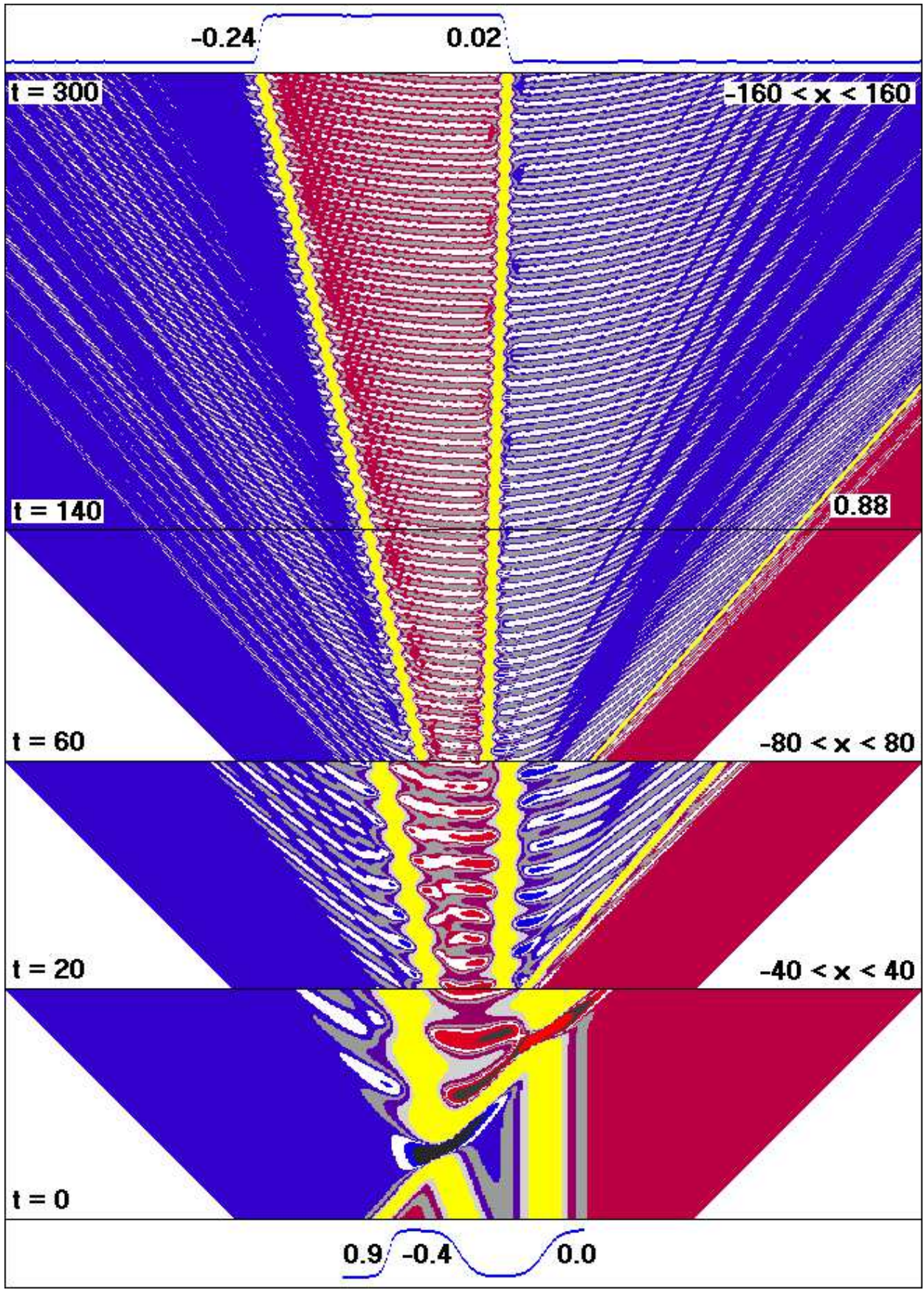}
\caption{Decay to three kinks}
\label{fig-4}
\end{center}
\end{figure}

\noindent
{\bf Numerical simulation.}
Our numerical experiments show a decay  of finite energy solutions to a finite set
of kinks and dispersive waves that corresponds to the asymptotics of \eqref {attN}.
One of the experiments is shown on Fig. \ref {fig-4}: a finite energy solution to equation (\ref {GL}) decays into three kinks.
Here the vertical line is the time axis, and the horizontal line is the   space axis;
the spatial scale redoubles at $ t = 20$ and at $ t = 60$.
Red color corresponds to $ \psi> 1+ \varepsilon $ values, blue color to $ \psi <-1-\varepsilon $ values,
and the yellow one to intermediate values $ -1 + \varepsilon <\psi <1-\varepsilon $, where $ \varepsilon> 0 $ is sufficiently small.
Thus, the yellow stripes represent  the kinks, while the blue and red zones outside the yellow stripes are filled with dispersive waves.
For $ t = 0 $, the solution begins with a rather chaotic behavior, when there are no visible kinks.
After 20 seconds, three separate kinks appear, which subsequently move almost uniformly.

\noindent {\bf The Lorentz contraction.}
The left kink moves to the left at a low speed $ v_1 \approx 0.24$, the central kink moves with a small velocity $ v_2 \approx 0.02 $,
and the right kink moves very fast with the speed $ v_3 \approx  0.88$.
The Lorentz spatial contraction $ \sqrt {1-v_k^2} $ is clearly visible on this picture:
the central kink is wide, the left one is a bit narrower,
and the right one is very narrow.

\smallskip

\noindent
{\bf The Einstein time delay.}
Also, the Einstein time delay is also very
prominent. 
Namely, all three kinks oscillate due to  the presence of  nonzero eigenvalue in the linearized equation at the kink:
substituting $ \psi(x,t) = S(x) + \varepsilon \varphi(x,t) $ in (\ref {GL}), we get in the first order
of $\varepsilon$
the linearized equation
\begin{equation}\label{lineq}
\ddot \varphi(x,t) = \varphi ''(x,t) -2 \varphi(x,t) -V(x) \varphi(x,t),
\end{equation}
 where the potential $V(x) = 3S^2(x) -3 =-\frac {3} {\cosh^2 {x} / {\sqrt2}}$ decays exponentially for large $ | x | $.
It is of great convenience
that for this potential the spectrum of the corresponding
\textit {Schr\"odinger operator} $ H: =-\frac {d^2} {dx^2} + 2 + V(x) $ is well known  \cite{Lamb80}.
Namely, the operator $ H $ is non-negative, and its continuous spectrum
is given by $ [2, \infty) $.
It turns out that $ H $ still has a two-point discrete spectrum: the points $ \lambda = 0 $ and $ \lambda = \frac32 $.
Exactly this nonzero eigenvalue  is responsible  for the pulsations that we observe for the central slow kink, with  the frequency
$ \omega_2 \approx \sqrt {\frac32} $ and period $ T_2 \approx 2 \pi\sqrt {\frac23}$. 
On the other hand, for  fast kinks, the ripples are much slower, i.e., the corresponding period is longer. 
This  time delay agrees numerically with the Lorentz formulas. These agreements
qualitatively confirm the relevance
of our numerical simulation. 

\smallskip

\noindent
{\bf Dispersive waves.}
The analysis of dispersive waves provides additional confirmation. Namely, the space outside the kinks on Fig. \ref {fig-4} is 
filled with dispersive waves, whose values are very close to $ \pm 1 $, with
accuracy $ 0.01 $.
These waves satisfy with high accuracy the linear Klein--Gordon equation which is obtained by linearization
of the Ginzburg--Landau equation (\ref {GL}) at the stationary solutions $ \psi_\pm \equiv \pm 1 $:
$$
\ddot \varphi(x,t) = \varphi ''(x,t) -2 \varphi(x,t).
$$
The corresponding dispersion relation $ \omega^2 = k^2 + 2 $ 
determines the group velocities of high-frequency wave packets:
\begin {equation} \label {ev}
\omega '(k)= \frac k {\sqrt {k^2 +2}} = \pm \frac {\sqrt {\omega^2 -2}} \omega.
\end {equation}
These wave packets are clearly visible on Fig. \ref {fig-4}
as straight lines whose propagation speeds converge to $ \pm 1 $.
This convergence is explained by the high-frequency limit  $ \omega '(k)\to \pm1 $  as $ \omega \to \pm \infty $.
For example, for dispersive waves emitted by central kink, the frequencies $ \omega = \pm n \omega_2 \to \pm \infty $
are generated by the polynomial nonlinearity in (\ref {GL}), see \cite{KK2020} for details.


The nonlinearity in (\ref {GL}) is chosen exactly because  of well-known spectrum of the linearized equation (\ref{lineq}). 
In numerical experiments \cite{KMV2004} we have also considered more general nonlinearities of the Ginzburg--Landau type. 
The results were qualitatively the same: for ``any'' initial data, the solution decays for large times into a sum of kinks and dispersive waves.
Numerically, this is clearly visible, but rigorous justification remains an open problem.

\subsection{Numerical observation of soliton asymptotics}
Besides the kinks,
 our numerical experiments \cite{KMV2004} have also resulted in the soliton-type asymptotics (\ref{attN}) 
and adiabatic effective dynamics  for solutions to the 1D relativistically-invariant nonlinear wave equations (\ref{GL}) with $F(\psi)=-U'(\psi)$ for the polynomial potentials
\begin{equation}\label{pop2}
  U(\psi)=a|\psi|^{2m}-b|\psi|^{2n},
\end{equation}
where $a,\,b>0$ and $m>n=2,\,3,\,\dotsc$. Respectively,
\begin{equation}\label{pop2f}
   F(\psi)=-2am|\psi|^{2m-2}\psi+2bn|\psi|^{2n-2}\psi.
\end{equation}
The parameters $a,\,m,\,b,\,n$ were taken as follows:
\begin{center}
\begin{tabular}{cccc}
$a$&$m$&$b$&$n$\\ \hline
1&3&0.61&2\\
10&4&2.1&2\\
10&6&8.75&5
\end{tabular}
\end{center}
We have considered various ``smooth''  initial functions $ \psi (x, 0), \dot\psi (x, 0) $ with  supports on the interval $[-20,20]$.
The second order finite-difference scheme with
$ \Delta x,\sim 0.01$, $ \Delta t \sim 0.001$
was employed. 
In all cases we have observed the asymptotics of type (\ref{attN}) with the numbers of solitons $0,1,3,5$ for $t> 100$.

\subsection{Adiabatic effective dynamics of relativistic solitons}

\begin{figure}[htbp]
\begin{center}
\includegraphics[width=1.0\columnwidth]{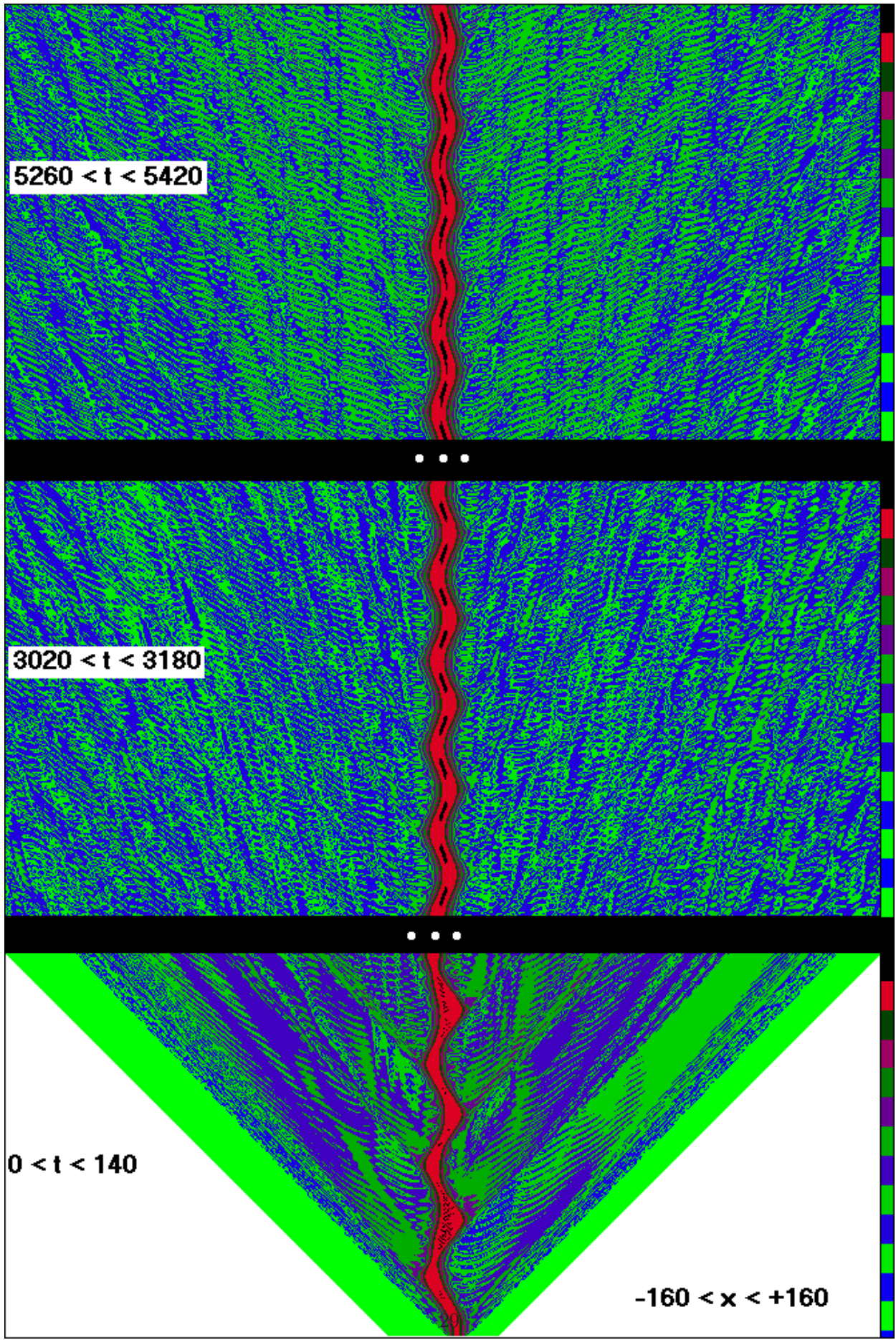}
\caption{Adiabatic effective dynamics of relativistic solitons}
\label{fig-7}
\end{center}
\end{figure}

In the numerical experiments \cite{KMV2004}
we also observed the adiabatic effective dynamics  for soliton-like solutions of the 1D equations with a~slowly varying external potential $V(x)$:
\begin{equation}\label{EQp}
  \ddot \psi(x,t)=\psi''(x,t)-\psi(x,t)+F(\psi(x,t))-V(x)\psi(x,t), \qquad x\in\R.
\end{equation}
This equation is formally equivalent to the Hamiltonian system with the Hamiltonian functional
\begin{equation}\label{HAMp}
  \mathcal{H}_V(\psi,\pi)=\int \Big[\frac12 |\pi(x)|^2+ \frac12 |\psi'(x)|^2+U(\psi(x))+\frac12 V(x) |\psi(x)|^2 \Big]\,dx.
\end{equation}
The soliton-like solutions
are
of the form 
\begin{equation}\label{swei}
  \psi(x,t)\approx e^{i\Theta(t)}\phi_{\omega(t)}(\gamma_{v(t)}(x-q(t))).
\end{equation}
Let us describe our numerical experiments which qualitatively confirm the adiabatic effective Hamiltonian dynamics
for the parameters $\Theta,\,\omega,\,q$, and $v$,
although its rigorous justification is still missing.
The effective dynamics of such type is proved in \cite{BCFJS}--\cite{S2010} 
for several Hamiltonian
models of PDEs coupled to a particle.

Figure~\ref{fig-7} represents solutions to equation (\ref{EQp}) with the potential (\ref{pop2}) with $a=10$, $m=6$, $b=8.75$, and $n = 5$.
We choose
\begin{equation}\label{v-cos}
V(x)=-0.2\cos (0.31 x)
\end{equation}
and the following initial conditions:
\begin{equation}\label{inco}
   \psi(x,0)= \phi_{\omega_0}(\gamma_{v_0}(x-q_0)), \qquad \dot \psi(x,0)=0,
\end{equation}
where $v_0=0$, $\omega_0=0.6$, and $q_0=5.0$.
We note that the initial state does not belong to  solitary manifold.
The effective width (half-amplitude) of the solitons is in the range $[4.4, 5.6 ]$.
It is quite small compared to the spatial period of the potential $2\pi/0.31 \sim 20$,
which is confirmed by numerical simulations shown on Figure~\ref{fig-7}. Namely,
\smallskip\\
$\bullet$ Blue and green colors represent a dispersive wave with values $|\psi(x,t)| <0.01$, while  red color represents a soliton with values $|\psi(x,t)|\in [0.4, 0.8]$.
\smallskip\\
$\bullet$ The soliton trajectory (a thick red meandering curve)
corresponds to oscillations of a classical particle in the potential $V(x)$.
\smallskip\\
$\bullet$ For $0<t<140$, the solution is far from the solitary manifold;
the radiation is intense.
\smallskip\\
$\bullet$ For $3020<t<3180 $, the solution approaches the solitary manifold; the radiation
weakens. The oscillation amplitude of the soliton is almost unchanged for a~long time, confirming  Hamilton-type dynamics.
\smallskip\\
$\bullet$ However, for $5260<t<5420$, the amplitude of the soliton oscillation is halved.
This suggests that at a~large time scale the deviation from the Hamiltonian
effective dynamics becomes essential.
Consequently, the effective dynamics gives a~good approximation only on the adiabatic time scale $t\sim \varepsilon^{-1}$.
\smallskip\\
$\bullet$ The deviation from the Hamiltonian dynamics is due to radiation, which plays the role of dissipation.
\smallskip\\
$\bullet$
The radiation is realized as  dispersive waves bringing the energy to the infinity.
The dispersive waves combine into uniformly moving bunches with discrete set of group velocities, as on Fig.~\ref{fig-4}.
The magnitude of solutions is of order $\sim 1$ on the trajectory of the soliton, while the values of the dispersive waves
is less than $0.01 $ for $t>200$, so that their energy density does not exceed $0.0001$.
The amplitude of the dispersive waves decays at large times.
In the limit $t\to\pm\infty$, the soliton converges to
a limit position
which corresponds to a local minimum
of the potential \eqref{v-cos}.

\end{document}